\documentclass[a4paper,10pt]{amsart}
\usepackage[german,english]{babel}
\usepackage[latin1]{inputenc}
\usepackage{ucs}
\usepackage{amsmath}
\usepackage{amsfonts}
\usepackage{amssymb}
\usepackage{amsthm}
\usepackage{fontenc}
\usepackage{graphicx}
\usepackage{url}
\usepackage{color}
\usepackage[T1]{fontenc}
\usepackage{a4wide}
\newtheorem{theorem}{Theorem}

\newtheorem{lemma}{Lemma}
\newtheorem{proposition}{Proposition}

\newtheorem{remark}{Remark}

\title[Functional a posteriori 
estimates for elliptic optimal control problems]{
A note on functional a posteriori
estimates for elliptic optimal control problems
}

\author{M. Wolfmayr}

\address[M. Wolfmayr]{
Johann Radon Institute for Computational and Applied Mathematics,
Austrian Academy of Sciences, Altenbergerstra\ss e 69, 4040 Linz, Austria}
\email{monika.wolfmayr@ricam.oeaw.ac.at}

\begin{document}

\begin{abstract}
In this work, 
new theoretical results on
functional type a posteriori estimates for elliptic optimal control problems
with control constraints are presented.
More precisely, we derive new, 
sharp, guaranteed and fully computable lower bounds for the cost functional
in addition to the already existing upper bounds.
Using both, the lower and the upper bounds,
we arrive at two-sided estimates for the cost functional.
We prove that these bounds 
finally lead to sharp,
guaranteed and fully computable upper estimates 
for the discretization error in the state and the control 
of the optimal control problem.
First numerical tests are presented 
confirming the efficiency of the a posteriori estimates derived.
\end{abstract}

\maketitle

%%%%%%%%%%%%%%%%%%%%%%%%%%%%%%%%%%%
\section{Introduction}
\label{Sec1:Intro}
%%%%%%%%%%%%%%%%%%%%%%%%%%%%%%%%%%%

During the last couple of decades, the optimization of systems
governed by partial differential equations (PDEs) 
has become more and more important in research and application, 
for which Lions has definitely paved the way with his work
\cite{WM:Lions:1971} in 1971.
Books considering
PDE-constrained optimization are, 
for instance,
Hinze et al. \cite{WM:HinzePinnauUlbrichUlbrich:2009},
Tr\"{o}ltzsch \cite{WM:Troeltzsch:2010},
Borz\`{i} and Schulz \cite{WM:BorziSchulz:2012}
and Leugering et al. \cite{WM:LeugeringEngellGriewankHinzeRannacherSchulzUlbrichUlbrich:2012,
WM:LeugeringBennerEngellGriewankHarbrechtHinzeRannacherUlbrich:2015}.
Besides the PDE-constraints, the optimization problems often include
control constraints given through 
a non-empty, convex and closed subset
of a Hilbert space.
In many cases, the set of admissible controls
is represented in terms of inequality constraints (or 
box constraints)
imposed on the controls.

A posteriori error estimates
and adaptive methods 
for elliptic optimal control problems
are the topic of many works, 
see, e.g.,
\cite{WM:BeckerKappRannacher:2000, WM:BeckerRannacher:2001,
WM:LiuYan:2001,
WM:BeckerBraackMeidnerRannacherVexler:2007,
WM:GuentherHinze:2008, WM:HintermuellerHoppe:2008,
WM:HintermuellerHoppeIliashKieweg:2008,
WM:VexlerWollner:2008, WM:HoppeKieweg:2010},
which are mainly on residual-type a posteriori error estimates.
Regarding a posteriori error estimates for control constrained
optimal control problems, we also refer to the recent works
\cite{WM:KohlsRoeschSiebert:2012, WM:KohlsRoeschSiebert:2014}.
There are different approaches to a posteriori error estimation.
Besides the residual-type estimates, there is also 
the class of functional type a posteriori error estimates.
These techniques were 
introduced by S.~Repin in the 90's, see, e.g.,
\cite{WM:Repin:1997, WM:Repin:1997b, WM:Repin:1999, WM:Repin:2000, WM:Repin:2000b, WM:NeittaanmaekiRepin:2004}.
Later the 
functional type a posteriori estimates were also considered and obtained for
optimal control problems,
see
\cite{WM:GaevskayaHoppeRepin:2006,
WM:GaevskayaHoppeRepin:2007}
as well as the books
\cite{WM:Repin:2008,
WM:MaliNeittaanmaekiRepin:2014}
and the references therein.
The benefit of functional type a posteriori estimates is that
- as the name indicates - they are only derived by functional methods. 
Therefore, these estimates 
do not depend on the mesh
and provide guaranteed upper bounds for the discretization errors.

In \cite{WM:GaevskayaHoppeRepin:2006},
sharp, guaranteed and fully computable upper bounds (majorants)
for cost functionals of distributed elliptic optimal control problems
were already presented.
Now, we want to complete the functional type a posteriori error analysis of
distributed elliptic optimal control problems
by the derivation of 
guaranteed and fully computable lower bounds (minorants) for the cost functionals.
The optimization problems of this work 
also include control constraints.
The presented minorants 
are not only fully computable but sharp.
Moreover, we prove that the new minorants and 
the already obtained majorants of \cite{WM:GaevskayaHoppeRepin:2006}
can be used in order to derive 
functional type a posteriori error
estimates for the discretization error in the state and the control.
The properties of the majorants and minorants are transferred to
the majorants for the discretization error, i.e., they are guaranteed,
fully computable and sharp.

The paper is organized as follows:
In Section~\ref{Sec2:EllOCP}, we discuss
the optimal control problem and its optimality system as well as present
the already known results on majorants 
for the cost functional.
Then, Section~\ref{Sec3:MinorantsForCostFunc} is devoted to the derivation of new,
fully computable minorants for the cost functional.
Both results, the older ones on majorants as well as the new ones on minorants,
together lead to 
guaranteed and fully computable upper bounds for the discretization error  
in the state and the control. These majorants are 
presented in Section~\ref{Sec4:FunctionalAPostErrorEstimates}.
They attain their exact lower bound on the exact solution of the optimal control problem
and, hence, are sharp.
In Section~\ref{Sec5:UnconstrainedCase}, we present the corresponding results for the optimal control problem in the case without any inequality constraints imposed on the control.
Section~\ref{Sec6:FEDiscretization} is devoted to the finite element discretization
of the optimal control problem and its iterative solution method 
in order to derive an approximation of the solution.
We present first numerical results in Section~\ref{Sec7:NumericalResults}, and 
finally draw some conclusions in Section~\ref{Sec8:Conclusions}.

%%%%%%%%%%%%%%%%%%%%%%%%%%%%%
\section{The Elliptic Optimal Control Problem}
\label{Sec2:EllOCP}
%%%%%%%%%%%%%%%%%%%%%%%%%%%%%

Let $\Omega \subset \mathbb{R}^d$, $d \in \{1,2,3\}$, be
a bounded Lipschitz domain with boundary $\Gamma := \partial \Omega$,
and, let us denote the state of our optimal control problem by $y$ and the control by $u$.
Given $y_d \in L^2(\Omega)$, $u_d \in L^2(\Omega)$, $f \in L^2(\Omega)$ and
$\lambda \in \mathbb{R}_+$,
we consider the following distributed elliptic optimal control problem: Minimize the cost functional
\begin{align}
\label{equation:minfunc:OCP}
 \mathcal{J}(y(v),v)
	:=  \frac{1}{2} \|y-y_d\|^2 
	+ \frac{\lambda}{2} \|v-u_d\|^2 
\end{align}
over $y \in V := H^1_0(\Omega)$ and $v \in U_{\text{ad}} \subset L^2(\Omega)$
subject to the boundary value problem
\begin{align}
\label{equation:forwardpde:OCP}
\begin{aligned}
 - \text{div} (\nu(\boldsymbol{x}) \nabla y(\boldsymbol{x})) 
 &= f(\boldsymbol{x}) + v(\boldsymbol{x}) \qquad &\boldsymbol{x} \in 
 \Omega, \\
 y(\boldsymbol{x}) &= 0 \qquad &\boldsymbol{x} \in 
 \Gamma. 
\end{aligned}
\end{align}
The set of admissible controls $U_{\text{ad}}$ 
is given by
\begin{align}
 \label{equation:boxConstr:OCP}
 U_{\text{ad}} = \{v \in L^2(\Omega): u_a \leq v \leq u_b 
 \, \text{ a.e. in } \Omega \},
\end{align}
where $u_a, u_b \in L^2(\Omega)$ 
and $u_a(\boldsymbol{x}) \leq u_b(\boldsymbol{x})$
for almost all $\boldsymbol{x} \in \Omega$.
The diffusion coefficient $\nu(\cdot)$ is assumed to be measurable, uniformly positive
and bounded, i.e., satisfies the assumptions
\begin{align}
 \label{assumptions:Nu}
 0 < \underline{\nu} \leq \nu(\boldsymbol{x}) \leq \overline{\nu}, \qquad 
 \boldsymbol{x} \in \Omega, 
\end{align}
where $\underline{\nu}$ and $\overline{\nu}$ are constants.
The positive regularization parameter $\lambda$ provides a weighting of the cost of the control 
in the cost functional $\mathcal{J}(\cdot,\cdot)$.
In these work, we denote by $(\cdot, \cdot)$ and $\| \cdot \|$ the inner products
and norms in $L^2(\Omega)$, respectively,
whereas the standard inner products and norms in $H^1(\Omega)$ are denoted by 
$(\cdot, \cdot)_1$ and $\| \cdot \|_1$, respectively.

Our goal is to approach to the desired state function $y_d$ and to the desired control $u_d$
as close as possible by finding a suitable control function $u$.
Note that, in many problems, the desired control is given by $u_d = 0$.
The optimal control problem (\ref{equation:minfunc:OCP})-(\ref{equation:boxConstr:OCP})
has a unique solution (for the proof, see, e.g.,
\cite{WM:HinzePinnauUlbrichUlbrich:2009, WM:Troeltzsch:2010})
that can be also derived via the optimality conditions.
Hence, the optimal solution of the optimality system
is equivalent to the solution of the optimal control problem
(\ref{equation:minfunc:OCP})-(\ref{equation:boxConstr:OCP}).
In order to formulate the optimality system, we consider
the following Lagrange functional for the minimization problem:
\begin{align}
 \label{equation:LagrangeFunctional}
 \mathcal{L}(y(v),v,p(v)):= \mathcal{J}(y(v),v) + \int_{\Omega} \big(
    \text{div}(\nu(\boldsymbol{x}) \nabla y(\boldsymbol{x})) 
    + f(\boldsymbol{x}) + v(\boldsymbol{x})\big) p(\boldsymbol{x}) \, d\boldsymbol{x},
\end{align}
where $p$ denotes the Lagrange multiplier (adjoint state).
The Lagrange functional (\ref{equation:LagrangeFunctional}) 
has to be understood in the weak sense. 
It has a saddle point, see, e.g., \cite{WM:HinzePinnauUlbrichUlbrich:2009, WM:Troeltzsch:2010}.
Hence, the corresponding (optimal) solution $(y,u,p) \in V \times U_{\text{ad}} \times V$
satisfies the system of (first order) necessary optimality conditions
  \begin{align*}
   \mathcal{L}_p(y(u), u, p(u)) &= 0, \\
   \mathcal{L}_y(y(u), u, p(u)) &= 0, \\
   (\mathcal{L}_u(y(u), u, p(u)), w - u)_{L^2(\Omega)} 
   &\geq 0
   \qquad \forall \, w \in U_{\text{ad}},
  \end{align*}
which can be written in the weak form as follows:
\begin{align}
\label{equation:KKTSysVF:1}
   \int_\Omega \nu \nabla y \cdot \nabla w \, d\boldsymbol{x}
   &= \int_\Omega (f + u) w \, d\boldsymbol{x}
   \qquad\, \forall \, w \in V, \\
\label{equation:KKTSysVF:2}
   \int_\Omega \nu \nabla p \cdot  \nabla w \, d\boldsymbol{x} &= \int_\Omega (y - y_d) w \, d\boldsymbol{x}
   \qquad \forall \, w \in V, \\
\label{equation:KKTSysVF:3}
   \int_\Omega (p + \lambda(u - u_d)) (w-u) \, d\boldsymbol{x} &\geq 0
   \qquad\qquad\qquad\qquad\, \forall \, w \in U_{\text{ad}},
\end{align}
for (optimal)
$y,p \in V$ and $u \in U_{\text{ad}}$.
\begin{remark}
 In the unconstrained case, i.e., $U_{\emph{ad}} = L^2(\Omega)$, the variational inequality
 (\ref{equation:KKTSysVF:3})
 of the optimality system simplifies to the equation $\mathcal{L}_u (y(u),u,p(u)) = 0$.
 Hence, the optimality condition (\ref{equation:KKTSysVF:3}) is then simplified to 
\begin{align}
 \label{equation:KKTSysVF:3Unconstrained}
    p + \lambda(u - u_d) = 0 \qquad\qquad\, \text{in } \, \Omega.
\end{align}
\end{remark}

\subsection*{Majorants for the cost functional}
Guaranteed and fully computable upper bounds
for the cost functional $\mathcal{J}$ have 
already been presented in \cite{WM:GaevskayaHoppeRepin:2006}, i.e.,
\begin{align}
\label{estimate:NotExactWithMajorant}
 \mathcal{J}(y(v),v) \leq \mathcal{J}^\oplus(\alpha,\beta;\eta,\boldsymbol{\tau},v) 
 \qquad \forall \, v \in U_{\text{ad}}
\end{align}
and for arbitrary
$\alpha, \beta > 0$, 
$\eta \in V$ and 
\begin{align*}
 \boldsymbol{\tau} \in H(\text{div},\Omega) :=
 \{\boldsymbol{\tau} \in [L^2(\Omega)]^d : \text{div} \, \boldsymbol{\tau} \in L^2(\Omega) \}.
\end{align*}
Although the admissible set is different in \cite{WM:GaevskayaHoppeRepin:2006}
and the elliptic boundary value problem is stated without a diffusion parameter $\nu(\cdot)$,
the majorants can be analogously derived for
the optimal control problem (\ref{equation:minfunc:OCP})-(\ref{equation:boxConstr:OCP}).
Including a diffusion parameter that meets the assumptions
(\ref{assumptions:Nu}) in the model problem (\ref{equation:forwardpde:OCP}),
the majorant $\mathcal{J}^\oplus$ is given by
\begin{align}
 \label{definition:majorant}
 \begin{aligned}
 \mathcal{J}^\oplus(\alpha,\beta;\eta,\boldsymbol{\tau},v):=
 &\frac{1+\alpha}{2}\|\eta - y_d\|^2
 + \frac{(1+\alpha)(1+\beta)C_F^2}{2\alpha\underline{\nu}^2} \|\boldsymbol{\tau}- \nu \nabla \eta\|^2 \\
 &+ \frac{(1+\alpha)(1+\beta)C_F^4}{2\alpha\beta\underline{\nu}^2}
 \|f+ v + \text{div} \, \boldsymbol{\tau}\|^2
 + \frac{\lambda}{2} \|v-u_d\|^2,
 \end{aligned}
\end{align}
where $C_F > 0$ is the constant
coming from the Friedrichs inequality.
The parameters 
$\alpha, \beta > 0$ have been introduced in order
to obtain a quadratic functional by applying Young's inequality.
The arbitrary functions $\eta \in V$ and $v \in U_{\text{ad}}$ 
can be taken as the approximate solutions of the optimal control problem
(\ref{equation:minfunc:OCP})-(\ref{equation:boxConstr:OCP})
and $\boldsymbol{\tau} \in H(\text{div},\Omega)$
represents the image of the
exact flux $\nu \nabla \eta$.
For the derivation of (\ref{definition:majorant}), the
following estimate for the approximation error has been used:
\begin{align}
 \label{definition:majorantBVP}
 \|\nabla y(v) - \nabla \eta\| \leq \frac{1}{\underline{\nu}} \left(\|\boldsymbol{\tau}-\nu \nabla \eta\|
 + C_F \|f+ v + \text{div} \, \boldsymbol{\tau}\|\right).
\end{align}
The majorant (\ref{definition:majorant}) provides a sharp upper bound of the cost functional,
if it is minimized over $\eta, \boldsymbol{\tau}, v$ and $\alpha, \beta > 0$, i.e.,
\begin{align}
\label{definition:majorantInf}
   \inf_{\substack{\eta \in V,\boldsymbol{\tau} \in H(\text{div},\Omega), \\
   v \in U_{\text{ad}},
   \alpha, \beta > 0}}
   \mathcal{J}^\oplus(\alpha,\beta;\eta,\boldsymbol{\tau},v) = \mathcal{J}(y(u),u),
\end{align}
since the infimum is attained for the optimal control $u$, its corresponding state $y(u)$
and its exact flux $\nu \nabla y(u)$, and for $\alpha$ going to zero.
Hence, (\ref{definition:majorantInf})  states that the exact lower bound of the
majorant (\ref{definition:majorant}) coincides with the optimal value of
the cost functional of the optimal control problem.
Therefore, we have the estimate
\begin{align}
\label{estimate:exactWithMajorant}
  \begin{aligned}
   \mathcal{J}(y(u),u)
   \leq \mathcal{J}^\oplus(\alpha,\beta;\eta,\boldsymbol{\tau},v) \qquad \,
   \forall \, \eta \in V, \,
   \boldsymbol{\tau} \in H(\text{div},\Omega), \,
   v \in U_{\text{ad}}, \, \alpha, \beta > 0,
  \end{aligned}
\end{align}
see \cite{WM:GaevskayaHoppeRepin:2006, WM:Repin:2008}.

%%%%%%%%%%%%%%%%%%%%%%%%%%%%
\section{Minorants for the Cost Functional} 
\label{Sec3:MinorantsForCostFunc}
%%%%%%%%%%%%%%%%%%%%%%%%%%%%

In this work, we enrich 
the derivation of guaranteed upper bounds for the discretization error
in the state and the control of 
problem (\ref{equation:minfunc:OCP})-(\ref{equation:boxConstr:OCP})
by obtaining fully computable lower bounds (minorants) for the cost functional
$\mathcal{J}$.
For any $\eta \in V$, we have that
\begin{align*}
 \mathcal{J}(y(v),v)
 &= \frac{1}{2} \|y-\eta\|^2 + \int_\Omega\left(y-\eta\right)\left(\eta-y_d\right)d\boldsymbol{x}
 + \frac{1}{2} \|\eta-y_d\|^2
 + \frac{\lambda}{2} \|v-u_d\|^2
\end{align*}
for all $v \in U_{\text{ad}}$.
Since $\frac{1}{2} \|y-\eta\|^2 \geq 0$, we can estimate $\mathcal{J}$ from below by
\begin{align}
\label{equation:costfunctionalExpandedEstBelow}
 \mathcal{J}(y(v),v)
  \geq \frac{1}{2} \|\eta-y_d\|^2
 + \frac{\lambda}{2} \|v-u_d\|^2 + \int_\Omega\left(y-\eta\right)\left(\eta-y_d\right)d\boldsymbol{x}.
\end{align}
Let $p_\eta \in V$ be the adjoint state corresponding to $\eta \in V$.
Hence, $p_\eta$ solves the equation
\begin{align}
 \label{equation:NOCforZeta}
   \int_\Omega \nu \nabla {p_\eta}\cdot  \nabla w \, d\boldsymbol{x} &= \int_\Omega (\eta - y_d) w \, d\boldsymbol{x}
   \qquad \forall \, w \in V.
\end{align}
By using (\ref{equation:NOCforZeta}), it follows for (\ref{equation:costfunctionalExpandedEstBelow}) that
\begin{align*}
 \mathcal{J}(y(v),v) &\geq 
 \frac{1}{2} \|\eta-y_d\|^2
 + \frac{\lambda}{2} \|v-u_d\|^2 +
 \int_\Omega\left(\nabla y - \nabla \eta\right) \cdot \nu \nabla {p_\eta}\, d\boldsymbol{x} \\
 &= \frac{1}{2} \|\eta-y_d\|^2
 + \frac{\lambda}{2} \|v-u_d\|^2 +
 \int_\Omega\left(\nu \nabla y - \nu \nabla \eta\right) \cdot \nabla {p_\eta}\, d\boldsymbol{x}. 
\end{align*}
Since $y = y(v)$ solves the variational formulation 
\begin{align}
  \label{equation:forwardpde:VF}
 \int_\Omega \nu \nabla y \cdot \nabla w \, d\boldsymbol{x}
 = \int_\Omega (f+v) w \, d\boldsymbol{x} \qquad \forall \, w \in V
\end{align}
of the boundary value problem (\ref{equation:forwardpde:OCP}),
we obtain
\begin{align*}
 \mathcal{J}(y(v),v)
 &\geq \frac{1}{2} \|\eta-y_d\|^2 + \frac{\lambda}{2} \|v-u_d\|^2
 + \int_\Omega (f+v) \, {p_\eta}\, d\boldsymbol{x}
 - \int_\Omega \nu \nabla \eta \cdot \nabla {p_\eta}\, d\boldsymbol{x}.
\end{align*}
For any $\boldsymbol{\tau} \in H(\text{div},\Omega)$, the identity
\begin{align}
\label{identity:Hdiv}
 \int_\Omega \text{div} \, \boldsymbol{\tau} \, w \, d\boldsymbol{x} 
 = -  \int_\Omega \boldsymbol{\tau} \cdot \nabla w \, d\boldsymbol{x} \qquad \forall \, w \in V
\end{align}
is valid, which yields
\begin{align*}
 \mathcal{J}(y(v),v)
 \geq \,& \frac{1}{2} \|\eta-y_d\|^2
 + \frac{\lambda}{2} \|v-u_d\|^2
 + \int_\Omega (f+v) \, {p_\eta}\, d\boldsymbol{x}
 - \int_\Omega \nu \nabla \eta \cdot \nabla {p_\eta}\, d\boldsymbol{x} \\
 &+ \int_\Omega \text{div} \, \boldsymbol{\tau} \, {p_\eta}\, d\boldsymbol{x} 
  + \int_\Omega \boldsymbol{\tau} \cdot \nabla {p_\eta}\, d\boldsymbol{x} \\
  = \,& \frac{1}{2} \|\eta-y_d\|^2
 + \frac{\lambda}{2} \|v-u_d\|^2 
 + \int_\Omega \left(f+v + \text{div} \, \boldsymbol{\tau}\right) {p_\eta}\, d\boldsymbol{x}
 + \int_\Omega \left(\boldsymbol{\tau} - \nu \nabla \eta \right) \cdot \nabla {p_\eta}\, d\boldsymbol{x} 
\end{align*}
for all $v \in U_{\text{ad}}$. 
Moreover,
\begin{align*}
  \begin{aligned}
 \mathcal{J}(y(u),u) =  
 \inf_{v \in U_{\text{ad}}} \mathcal{J}(y(v),v)
 \geq &\,
  \frac{1}{2} \|\eta-y_d\|^2
  + \int_\Omega \left(f+ \text{div} \, \boldsymbol{\tau}\right) {p_\eta}\, d\boldsymbol{x}
  + \int_\Omega \left(\boldsymbol{\tau} - \nu \nabla \eta \right) \cdot \nabla {p_\eta}\, d\boldsymbol{x} \\
  &+ \inf_{v \in U_{\text{ad}}}\left(
  \int_\Omega v \, {p_\eta}\, d\boldsymbol{x} 
  + \frac{\lambda}{2} \|v-u_d\|^2 
  \right).
  \end{aligned}
\end{align*}
The control 
$v \in U_{\text{ad}}$ can be defined via the projection formula
\begin{align*}
 \mathbb{P}_{[a,b]}(u):=\min\{b,\max\{a,u\}\}
\end{align*}
for all $a, b \in \mathbb{R}$ with $a \leq b$ and $u \in \mathbb{R}$,
which projects $\mathbb{R}$ on the interval $[a,b]$.
It is given by
\begin{align*}
  v(\boldsymbol{x}) = \mathbb{P}_{[u_a(\boldsymbol{x}),u_b(\boldsymbol{x})]}
 \{u_d(\boldsymbol{x}) - \frac{1}{\lambda} p(\boldsymbol{x}) \},
\end{align*}
where $p = p(v)$ is here the adjoint state corresponding to $v$ and $y(v)$.
Let us denote by $v_{p_\eta} \in U_{\text{ad}}$ 
the control corresponding to the adjoint state ${p_\eta}$, i.e., 
\begin{align}
\label{definition:projection}
 v_{p_\eta}(\boldsymbol{x}) = \mathbb{P}_{[u_a(\boldsymbol{x}),u_b(\boldsymbol{x})]}
 \{u_d(\boldsymbol{x}) - \frac{1}{\lambda} {p_\eta}(\boldsymbol{x}) \}
\end{align}
for almost every $\boldsymbol{x} \in \Omega$. 
By adding and subtracting $v_{p_\eta}$, we obtain
\begin{align*}
 \mathcal{J}(y(v),v)
  \geq \,& \frac{1}{2} \|\eta-y_d\|^2
 + \frac{\lambda}{2} \|v_{p_\eta}-u_d\|^2 + \frac{\lambda}{2} \|v - v_{p_\eta}\|^2
 + \lambda \int_\Omega \left(v - v_{p_\eta}\right) \left(v_{p_\eta}-u_d\right) \, d\boldsymbol{x} \\
 &+ \int_\Omega \left(v- v_{p_\eta}\right) {p_\eta}\, d\boldsymbol{x}
 + \int_\Omega \left(f+ v_{p_\eta}+ \text{div} \, \boldsymbol{\tau}\right) {p_\eta}\, d\boldsymbol{x}
 + \int_\Omega \left(\boldsymbol{\tau} - \nu \nabla \eta \right) \cdot \nabla {p_\eta}\, d\boldsymbol{x} 
\end{align*}
for all $v \in U_{\text{ad}}$. Since $\frac{\lambda}{2} \|v - v_{p_\eta}\|^2 \geq 0$, we get that
\begin{align*}
 \mathcal{J}(y(v),v)
  \geq \,& \frac{1}{2} \|\eta-y_d\|^2
 + \frac{\lambda}{2} \|v_{p_\eta}-u_d\|^2 
 + \int_\Omega \left(v - v_{p_\eta}\right) \left({p_\eta}+ \lambda (v_{p_\eta}-u_d)\right) \, d\boldsymbol{x} \\
 &+ \int_\Omega \left(f+ v_{p_\eta}+ \text{div} \, \boldsymbol{\tau}\right) {p_\eta}\, d\boldsymbol{x}
 + \int_\Omega \left(\boldsymbol{\tau} - \nu \nabla \eta \right) \cdot \nabla {p_\eta}\, d\boldsymbol{x}.
\end{align*}
Due to the variational inequality 
\begin{align}
\label{inequality:VarIneqZeta}
 \int_\Omega \left(v - v_{p_\eta}\right) \left({p_\eta}
 + \lambda (v_{p_\eta}-u_d)\right) \, d\boldsymbol{x} \geq 0 \qquad
 \forall \, v \in U_{\text{ad}},
\end{align}
we obtain
\begin{align*}
 \mathcal{J}(y(v),v)
  \geq \,& \frac{1}{2} \|\eta-y_d\|^2
 + \frac{\lambda}{2} \|v_{p_\eta}-u_d\|^2
 + \int_\Omega \left(f+ v_{p_\eta}+ \text{div} \, \boldsymbol{\tau}\right) {p_\eta}\, d\boldsymbol{x}
 + \int_\Omega \left(\boldsymbol{\tau} - \nu \nabla \eta \right) \cdot \nabla {p_\eta}\, d\boldsymbol{x}.
\end{align*}
Now, let us introduce an arbitrary $\zeta \in V$ and its corresponding control $v_\zeta \in U_{\text{ad}}$,
which can be computed by the projection formula as follows:
\begin{align}
\label{definition:projectionDiscrete}
 v_{\zeta}(\boldsymbol{x}) = \mathbb{P}_{[u_a(\boldsymbol{x}),u_b(\boldsymbol{x})]}
 \{u_d(\boldsymbol{x}) - \frac{1}{\lambda} \zeta(\boldsymbol{x}) \}.
\end{align}
Hence, the
following variational inequality is satisfied: 
\begin{align}
\label{inequality:VarIneqZetah}
 \int_\Omega \left(v - v_{\zeta}\right) \left(\zeta + \lambda (v_{\zeta} -u_d)\right) \, d\boldsymbol{x} \geq 0 \qquad
 \forall \, v \in U_{\text{ad}}.
\end{align}
Adding and subtracting $v_{\zeta} \in U_{\text{ad}}$ leads to
\begin{align*}
 \mathcal{J}(y(v),v)
  \geq \,& \frac{1}{2} \|\eta-y_d\|^2
 + \frac{\lambda}{2} \|v_{\zeta} -u_d\|^2
 + \lambda \int_\Omega \left(v_{p_\eta}- v_{\zeta}\right) \left(v_{\zeta} -u_d\right) \, d\boldsymbol{x} \\
 &+ \int_\Omega \left(f+ v_{\zeta} + \text{div} \, \boldsymbol{\tau}\right) {p_\eta}\, d\boldsymbol{x}
 + \int_\Omega \left(\boldsymbol{\tau} - \nu \nabla \eta \right) \cdot \nabla {p_\eta}\, d\boldsymbol{x}
 + \int_\Omega \left(v_{p_\eta}-v_{\zeta}\right) {p_\eta}\, d\boldsymbol{x} \\
 \geq \,& \frac{1}{2} \|\eta-y_d\|^2
 + \frac{\lambda}{2} \|v_{\zeta} -u_d\|^2
 + \int_\Omega \left(v_{p_\eta}- v_{\zeta}\right)({p_\eta}-\zeta) \, d\boldsymbol{x} \\
 &+ \int_\Omega \left(f+ v_{\zeta} + \text{div} \, \boldsymbol{\tau}\right) {p_\eta}\, d\boldsymbol{x}
 + \int_\Omega \left(\boldsymbol{\tau} - \nu \nabla \eta \right) \cdot \nabla {p_\eta}\, d\boldsymbol{x},
\end{align*}
since $\frac{\lambda}{2} \|v_{p_\eta}- v_{\zeta}\|^2 \geq 0$ and the 
variational inequality (\ref{inequality:VarIneqZetah}) is valid. 
Next, we add and subtract the arbitrary 
$\zeta \in V$ leading to the following estimate:
\begin{align}
\label{inequality:estimateJWitherror}
\begin{aligned}
 \mathcal{J}(y(v),v)
 \geq \,& \frac{1}{2} \|\eta-y_d\|^2
 + \frac{\lambda}{2} \|v_{\zeta} -u_d\|^2
 + \int_\Omega \left(v_{p_\eta}- v_{\zeta}\right)({p_\eta}-\zeta) \, d\boldsymbol{x} \\
 &+ \int_\Omega \left(f+ v_{\zeta} + \text{div} \, \boldsymbol{\tau}\right) \zeta \, d\boldsymbol{x}
 + \int_\Omega \left(\boldsymbol{\tau} - \nu \nabla \eta \right) \cdot \nabla \zeta \, d\boldsymbol{x} \\
 &+ \int_\Omega \left(f+ v_{\zeta} + \text{div} \, \boldsymbol{\tau}\right) ({p_\eta}-\zeta) \, d\boldsymbol{x}
 + \int_\Omega \left(\boldsymbol{\tau} - \nu \nabla \eta \right) \cdot \nabla ({p_\eta}-\zeta) \, d\boldsymbol{x}.
\end{aligned}
\end{align}
The following result provides an estimate for the error in the control:
\begin{lemma}
\label{lemma:estimate:control}
 Let $v_{p_\eta}\in U_{\emph{ad}}$ and $v_{\zeta} \in U_{\emph{ad}}$ satisfy the variational inequalities
 (\ref{inequality:VarIneqZeta}) and (\ref{inequality:VarIneqZetah}), respectively.
 Then, the error between $v_{p_\eta}$ and $v_{\zeta}$ can be estimated by
 \begin{align}
 \label{estimate:control}
  \|v_{p_\eta}- v_{\zeta}\| \leq \frac{1}{\lambda} \|{p_\eta}- \zeta\|.
 \end{align}
 \begin{proof}
  Adding the variational inequalities
  \begin{align*}
  \int_\Omega \left(v_{p_\eta}- v_{\zeta}\right) \left(\zeta + \lambda (v_{\zeta} -u_d)\right) \, d\boldsymbol{x} \geq 0
  \qquad \text{and} \qquad
  \int_\Omega \left(v_{\zeta} - v_{p_\eta}\right) \left({p_\eta}+ \lambda (v_{p_\eta}-u_d)\right) \, d\boldsymbol{x} \geq 0
  \end{align*}  
  yields the inequality
  \begin{align*}
  \int_\Omega \left(v_{p_\eta}- v_{\zeta}\right)
  \left(\zeta-{p_\eta}\right) \, d\boldsymbol{x} \geq \lambda
  \int_\Omega (v_{p_\eta} - v_{\zeta})^2 \, d\boldsymbol{x}.
  \end{align*}
  By applying the Cauchy-Schwarz inequality, we obtain
  \begin{align*}
  \|v_{p_\eta}- v_{\zeta}\|\|{p_\eta}- \zeta\| \geq \int_\Omega \left(v_{p_\eta}- v_{\zeta}\right)
  \left(\zeta-{p_\eta}\right) \, d\boldsymbol{x} \geq \lambda
  \|v_{p_\eta} - v_{\zeta}\|^2,
  \end{align*}
  which finally leads to the estimate (\ref{estimate:control}).
 \end{proof}
\end{lemma}
In order to formulate, a computable lower bound for the cost functional, we need to prove a computable
upper bound for the error in the adjoint state, which is presented in the following theorem:
\begin{theorem}
\label{theorem:estimate:adjointState}
 Let $y_d \in L^2(\Omega)$ be given and let ${p_\eta}\in V$ meet equation (\ref{equation:NOCforZeta}) 
 with $\eta \in V$ and $\nu(\cdot)$ satisfying assumption (\ref{assumptions:Nu}). For any 
 $\zeta\in V$, we have that
 \begin{align}
  \label{estimate:adjointState}
  \|\nabla ({p_\eta}- \zeta)\| \leq \frac{1}{\underline{\nu}}
  \left(C_F \, \|\eta - y_d + \emph{div} \, \boldsymbol{\rho}\| + \|\boldsymbol{\rho} - \nu \nabla \zeta\| \right),
 \end{align}
 where $\boldsymbol{\rho} \in H(\emph{div},\Omega)$
 and $C_F > 0$ is the constant coming from the Friedrichs inequality.
 \begin{proof}
 Since the bilinear form of problem (\ref{equation:NOCforZeta}) is elliptic with ellipticity constant
 $\underline{\nu}$ and applying the Cauchy-Schwarz and Friedrichs inequalities, we get that
  \begin{align*}
    \underline{\nu} \|\nabla ({p_\eta}- \zeta)\|
    &\leq \sup_{0 \not=w \in V} \frac{\int_\Omega \nu \nabla ({p_\eta}- \zeta) \cdot  \nabla w \, d\boldsymbol{x}}{\|\nabla w\|}
    = \sup_{0 \not=w \in V} \frac{\int_\Omega (\eta - y_d) w
    - \nu \nabla \zeta\cdot  \nabla w \, d\boldsymbol{x}}{\|\nabla w\|} \\
    &= \sup_{0 \not=w \in V} \frac{\int_\Omega (\eta - y_d + \text{div} \, \boldsymbol{\rho}) w
    +(\boldsymbol{\rho} - \nu \nabla \zeta) \cdot  \nabla w \, d\boldsymbol{x}}{\|\nabla w\|} \\
    &\leq \sup_{0 \not=w \in V} \frac{\|\eta - y_d + \text{div} \, \boldsymbol{\rho}\| \|w\|
    +\|\boldsymbol{\rho} - \nu \nabla \zeta\| \|\nabla w\|}{\|\nabla w\|} \\
    &\leq C_F \, \|\eta - y_d + \text{div} \, \boldsymbol{\rho}\|
    +\|\boldsymbol{\rho} - \nu \nabla \zeta\|
  \end{align*}
  where $\boldsymbol{\rho} \in H(\text{div},\Omega)$ satisfies identity (\ref{identity:Hdiv}).
  Hence, it follows 
  the estimate (\ref{estimate:adjointState}).
 \end{proof}
\end{theorem}
By using the Cauchy-Schwarz and Friedrichs inequalities as well as Lemma~\ref{lemma:estimate:control},
we can further estimate the inequality (\ref{inequality:estimateJWitherror}) from below as follows:
\begin{align*}
\begin{aligned}
 \mathcal{J}(y(v),v)
 \geq \,& \frac{1}{2} \|\eta-y_d\|^2
 + \frac{\lambda}{2} \|v_{\zeta} -u_d\|^2
 + \int_\Omega \left(f+ v_{\zeta} + \text{div} \, \boldsymbol{\tau}\right) {\zeta} \, d\boldsymbol{x}
 + \int_\Omega \left(\boldsymbol{\tau} - \nu \nabla \eta \right) \cdot \nabla {\zeta} \, d\boldsymbol{x} \\
 &- \|v_{p_\eta}- v_{\zeta}\|\|{p_\eta}-{\zeta}\|
 -\|f+ v_{\zeta} + \text{div} \, \boldsymbol{\tau}\|\|{p_\eta}-{\zeta}\|
 -\|\boldsymbol{\tau} - \nu \nabla \eta \| \|\nabla ({p_\eta}-{\zeta}) \| \\
 \geq \,& \frac{1}{2} \|\eta-y_d\|^2
 + \frac{\lambda}{2} \|v_{\zeta} -u_d\|^2
 + \int_\Omega \left(f+ v_{\zeta} + \text{div} \, \boldsymbol{\tau}\right) {\zeta} \, d\boldsymbol{x}
 + \int_\Omega \left(\boldsymbol{\tau} - \nu \nabla \eta \right) \cdot \nabla {\zeta} \, d\boldsymbol{x} \\
 &- \frac{1}{\lambda} \|{p_\eta}- {\zeta}\|^2
 -(C_F \, \|f+ v_{\zeta} + \text{div} \, \boldsymbol{\tau}\|
 +\|\boldsymbol{\tau} - \nu \nabla \eta \|) \|\nabla ({p_\eta}-{\zeta}) \| \\
 \geq \,& \frac{1}{2} \|\eta-y_d\|^2
 + \frac{\lambda}{2} \|v_{\zeta} -u_d\|^2
 + \int_\Omega \left(f+ v_{\zeta} + \text{div} \, \boldsymbol{\tau}\right) {\zeta} \, d\boldsymbol{x} 
 + \int_\Omega \left(\boldsymbol{\tau} - \nu \nabla \eta \right) \cdot \nabla {\zeta} \, d\boldsymbol{x} \\
 &- \frac{C_F^2}{\lambda} \|\nabla({p_\eta}- {\zeta})\|^2 
 -(C_F \, \|f+ v_{\zeta} + \text{div} \, \boldsymbol{\tau}\|
 +\|\boldsymbol{\tau} - \nu \nabla \eta \|) \|\nabla ({p_\eta}-{\zeta}) \|.
\end{aligned}
\end{align*}
Now, applying Theorem~\ref{theorem:estimate:adjointState} yields the following estimate:
\begin{align}
\label{inequality:finalminorant}
\begin{aligned}
 \mathcal{J}(y(v),v)
 \geq
 \mathcal{J}^\ominus(\eta,{\zeta},\boldsymbol{\tau},\boldsymbol{\rho},v_{\zeta}) \qquad
 \forall \eta, {\zeta} \in V, \boldsymbol{\tau},\boldsymbol{\rho} \in H(\text{div},\Omega), 
 v_{\zeta} = \mathbb{P}_{[u_a,u_b]} \{u_d - \frac{1}{\lambda} {\zeta} \}
\end{aligned}
\end{align}
with the minorant
\begin{align}
 \label{definition:minorant}
\begin{aligned}
 \mathcal{J}^\ominus(\eta,{\zeta}&,\boldsymbol{\tau},\boldsymbol{\rho},v_{\zeta})
 = \frac{1}{2} \|\eta-y_d\|^2
 + \frac{\lambda}{2} \|v_{\zeta} -u_d\|^2
 + \int_\Omega \left(f+ v_{\zeta} + \text{div} \, \boldsymbol{\tau}\right) {\zeta} \, d\boldsymbol{x} \\
 &+ \int_\Omega \left(\boldsymbol{\tau} - \nu \nabla \eta \right) \cdot \nabla {\zeta} \, d\boldsymbol{x} 
 - \frac{1}{\underline{\nu}} \left(C_F \, \|\eta - y_d + \text{div} \, \boldsymbol{\rho}\| + \|\boldsymbol{\rho} - \nu \nabla {\zeta}\| \right) \\
 &\times \left(\frac{C_F^3}{\lambda \underline{\nu}} \|\eta - y_d + \text{div} \, \boldsymbol{\rho}\|
 + \frac{C_F^2}{\lambda \underline{\nu}} \|\boldsymbol{\rho} - \nu \nabla {\zeta}\| 
 + C_F \|f+ v_{\zeta} + \text{div} \, \boldsymbol{\tau}\|
 + \|\boldsymbol{\tau} - \nu \nabla \eta \| \right),
\end{aligned}
\end{align}
where $\boldsymbol{\tau}, \boldsymbol{\rho} \in H(\text{div},\Omega)$.
Note that the minorant is fully computable.
\begin{theorem} 
 The exact upper bound of the minorant $\mathcal{J}^\ominus$
 defined in (\ref{definition:minorant})
 coincides with the optimal value of the cost functional of problem
 (\ref{equation:minfunc:OCP})-(\ref{equation:boxConstr:OCP}), or, equivalently, of
 the optimality system (\ref{equation:KKTSysVF:1})-(\ref{equation:KKTSysVF:3}),
 i.e.,
\begin{align}
 \label{definition:minorantSup}
  \sup_{\substack{\eta,{\zeta} \in V, \boldsymbol{\tau},\boldsymbol{\rho} \in H(\emph{div},\Omega), \\
    v_{\zeta} = \mathbb{P}_{[u_a,u_b]} \{u_d - \frac{1}{\lambda} {\zeta} \}}}
  \mathcal{J}^\ominus(\eta,{\zeta},\boldsymbol{\tau},\boldsymbol{\rho},v_{\zeta})
  = \mathcal{J}(y(u),u).
\end{align}
\begin{proof}
The estimate (\ref{inequality:finalminorant}) is valid for all $v \in U_{\text{ad}}$.
Hence, also for the exact solution $u$, i.e.,
\begin{align}
 \label{estimate:exactWithMinorant}
  \begin{aligned}
 \mathcal{J}(y(u),u) = 
 \inf_{v \in U_{\text{ad}}} \mathcal{J}(y(v),v)
 \geq \mathcal{J}^\ominus(\eta,{\zeta},\boldsymbol{\tau},\boldsymbol{\rho},v_{\zeta})
  \end{aligned}
\end{align}
for all $\eta,\zeta \in V$, $\boldsymbol{\tau},\boldsymbol{\rho} \in H(\text{div},\Omega)$ and for
the control $v_{\zeta} \in U_{\text{ad}}$, which depends on $\zeta$ and can be computed by
the projection formula (\ref{definition:projectionDiscrete}).
For the exact solution $v_{\zeta} = u$, 
$\eta = y(u)$, $\zeta= p(u)$, $\boldsymbol{\tau} = \nu \nabla y(u)$
and $\boldsymbol{\rho} = \nu \nabla p(u)$,
the estimate is sharp, i.e.,
\begin{align*}
 \mathcal{J}^\ominus(y(u),p(u),\nu \nabla y(u),\nu \nabla p(u),u)
 = &\, \frac{1}{2} \|y-y_d\|^2
 + \frac{\lambda}{2} \|u-u_d\|^2 = \mathcal{J}(y(u),u).
\end{align*}
\end{proof}
\end{theorem}

%%%%%%%%%%%%%%%%%%%%%%%%%%%%%%%%%%
\section{A Posteriori Error Estimates for Control and State}
\label{Sec4:FunctionalAPostErrorEstimates}
%%%%%%%%%%%%%%%%%%%%%%%%%%%%%%%%%%

In this section, we will 
derive guaranteed upper bounds for the discretization errors of the control 
and the state measured in the following norm:
\begin{align}
 \label{definition:combinedNorm}
 |||u-v|||^2 := \frac{1}{2} \|y(u) - y(v)\|^2 + \frac{\lambda}{2} \|u-v\|^2,
\end{align}
making use of the ideas 
based on the work by Mikhlin \cite{WM:Mikhlin:1964} but 
generalized for the class of optimal control problems, see also \cite{WM:Repin:2008}.
%--------------------------------------------------------------------------------------------------
\begin{theorem}
 \label{theorem:DiffCombNorm}
 For any control function $v \in U_{\text{ad}}$,
 we have the estimate 
 \begin{align}
  \label{equation:DiffCombNorm}
  |||u-v|||^2 \leq 
  \mathcal{J}(y(v),v) - \mathcal{J}(y(u),u).
 \end{align}
\begin{proof}
 We compute the difference
 \begin{align*}
  \mathcal{J}(y(v)&,v) - \mathcal{J}(y(u),u)
  = \frac{1}{2} \|y(v)-y_d\|^2 - \frac{1}{2} \|y(u)-y_d\|^2
    + \frac{\lambda}{2} \|v-u_d\|^2 - \frac{\lambda}{2} \|u-u_d\|^2 \\
  = &\,\frac{1}{2} \int_\Omega (y(v)+y(u)-2y_d)(y(v)-y(u)) \, d\boldsymbol{x} 
    + \frac{\lambda}{2} \int_\Omega (v+u-2u_d)(v-u) \, d\boldsymbol{x} \\
  = &\,\frac{1}{2} \int_\Omega (y(v)-y(u)+2 y(u)-2 y_d)(y(v)-y(u)) \, d\boldsymbol{x} 
    + \frac{\lambda}{2} \int_\Omega (v-u+2 u-2 u_d)(v-u) \, d\boldsymbol{x} \\
  = &\,\frac{1}{2} \|y(u) - y(v)\|^2 + \int_\Omega (y(u)-y_d)(y(v)-y(u)) \, d\boldsymbol{x} 
    + \frac{\lambda}{2} \|u-v\|^2 + \lambda \int_\Omega (u-u_d)(v-u) \, d\boldsymbol{x}.
 \end{align*}
 Since the exact adjoint state $p(u) \in V$ fulfills
 (\ref{equation:KKTSysVF:2}), we obtain the equation
 \begin{align*}
  \mathcal{J}(y(v&),v) - \mathcal{J}(y(u),u) \\
  = &\,\frac{1}{2} \|y(u) - y(v)\|^2 + \int_\Omega \nu \nabla p(u)(\nabla y(v)- \nabla y(u)) \, d\boldsymbol{x} 
    + \frac{\lambda}{2} \|u-v\|^2 + \lambda \int_\Omega (u-u_d)(v-u) \, d\boldsymbol{x} \\
  = &\,\frac{1}{2} \|y(u) - y(v)\|^2 + \int_\Omega \nabla p(u)(\nu \nabla y(v)- \nu \nabla y(u)) \, d\boldsymbol{x} 
    + \frac{\lambda}{2} \|u-v\|^2 + \lambda \int_\Omega (u-u_d)(v-u) \, d\boldsymbol{x}.
 \end{align*}
 From equations (\ref{equation:KKTSysVF:1}) and
 (\ref{equation:forwardpde:VF}) follows that
 \begin{align*}
  \mathcal{J}(y(v)&,v) - \mathcal{J}(y(u),u) \\
  = &\,\frac{1}{2} \|y(u) - y(v)\|^2 + \int_\Omega p(u)(f+v-f-u) \, d\boldsymbol{x} 
    + \frac{\lambda}{2} \|u-v\|^2 + \lambda \int_\Omega (u-u_d)(v-u) \, d\boldsymbol{x} \\
  = &\,|||u-v|||^2
    + \int_\Omega (p(u)+\lambda(u-u_d))(v-u) \, d\boldsymbol{x}.
 \end{align*}
 Since the variational inequality (\ref{equation:KKTSysVF:3}) is satisfied, i.e.,
 \begin{align*}
  \int_\Omega (p(u)+\lambda(u-u_d))(v-u) \, d\boldsymbol{x} \geq 0 \qquad \forall \, v \in U_{\text{ad}},
 \end{align*} 
 we finally obtain the estimate (\ref{equation:DiffCombNorm}).
\end{proof}
\end{theorem}
%--------------------------------------------------------------------------------------------------
\begin{theorem}
\label{theorem:majorantCombNorm}
For any $\zeta \in V$, let $v_{\zeta} \in U_{\text{ad}}$
be given by 
the projection formula (\ref{definition:projectionDiscrete}).
Then, we obtain the following error majorant:
\begin{align}
 \label{definition:majorantCombinedNorm}
  |||u-v_\zeta|||^2 \leq \mathcal{M}^\oplus(\alpha,\beta;\eta,\zeta,\boldsymbol{\tau},\boldsymbol{\rho},v_{\zeta})
  := \mathcal{J}^\oplus(\alpha,\beta;\eta,\boldsymbol{\tau},v_\zeta)
  - \mathcal{J}^\ominus(\eta,{\zeta},\boldsymbol{\tau},\boldsymbol{\rho},v_\zeta)
\end{align}
with
\begin{align*}
 \mathcal{M}^\oplus(\alpha,\beta;\eta,\zeta&,\boldsymbol{\tau},\boldsymbol{\rho},v_{\zeta})
 = \frac{\alpha}{2}\|\eta - y_d\|^2
 + \frac{(1+\alpha)(1+\beta)C_F^2}{2\alpha\underline{\nu}^2} \|\boldsymbol{\tau}- \nu \nabla \eta\|^2 \\
 &+ \frac{(1+\alpha)(1+\beta)C_F^4}{2\alpha\beta\underline{\nu}^2} \|f+ v_{\zeta} + \emph{div} \, \boldsymbol{\tau}\|^2
 - \int_\Omega \left(f+ v_{\zeta} + \emph{div} \, \boldsymbol{\tau}\right) {\zeta} \, d\boldsymbol{x} \\
 &- \int_\Omega \left(\boldsymbol{\tau} - \nu \nabla \eta \right) \cdot \nabla {\zeta} \, d\boldsymbol{x}
 + \frac{1}{\underline{\nu}} \left(C_F \, \|\eta - y_d + \emph{div} \, \boldsymbol{\rho}\|
 + \|\boldsymbol{\rho} - \nu \nabla {\zeta}\| \right) \\
 &\times \left(\frac{C_F^3}{\lambda \underline{\nu}} \|\eta - y_d + \emph{div} \, \boldsymbol{\rho}\|
 + \frac{C_F^2}{\lambda \underline{\nu}} \|\boldsymbol{\rho} - \nu \nabla {\zeta}\| 
 + C_F \|f+ v_{\zeta} + \emph{div} \, \boldsymbol{\tau}\|
 + \|\boldsymbol{\tau} - \nu \nabla \eta \| \right)
\end{align*}
for arbitrary
$\eta \in V$, 
$\boldsymbol{\tau}, \boldsymbol{\rho} \in H(\emph{div},\Omega)$
and 
$\alpha,\beta > 0$.
\begin{proof}
  Using (\ref{equation:DiffCombNorm}), we obtain the estimate
  \begin{align*}
   |||u-v_\zeta|||^2 \leq \mathcal{J}(y(v_\zeta),v_\zeta) - \mathcal{J}(y(u),u).
  \end{align*}
  Applying 
  (\ref{estimate:NotExactWithMajorant}) and (\ref{estimate:exactWithMinorant}) 
  finally leads to
  the estimate (\ref{definition:majorantCombinedNorm}).
\end{proof}
\end{theorem}
%--------------------------------------------------------------------------------------------------
\begin{proposition}
\label{proposition:overallmajorant}
 The majorant
 $\mathcal{M}^\oplus(\alpha,\beta;\eta,\zeta,\boldsymbol{\tau},\boldsymbol{\rho},v_{\zeta})$ defined in
 (\ref{definition:majorantCombinedNorm}) attains the exact lower bound on the exact solution of
 the optimal control problem (\ref{equation:minfunc:OCP})-(\ref{equation:boxConstr:OCP}), or,
 equivalently, of
 the optimality system (\ref{equation:KKTSysVF:1})-(\ref{equation:KKTSysVF:3}),
 i.e.,
\begin{align*}
    \inf_{\substack{\eta,{\zeta} \in V, \boldsymbol{\tau},\boldsymbol{\rho} \in H(\emph{div},\Omega), \\
    v_{\zeta} = \mathbb{P}_{[u_a,u_b]} \{u_d - \frac{1}{\lambda} {\zeta} \}, \alpha,\beta > 0}}
  \mathcal{M}^\oplus(\alpha,\beta;\eta,\zeta,\boldsymbol{\tau},\boldsymbol{\rho},v_{\zeta}) = 0.
\end{align*}
 The infimum is attained for $v_\zeta=u$, $\eta = y(u)$, ${\zeta}= p(u)$,
 $\boldsymbol{\tau} = \nu \nabla y(u)$ and $\boldsymbol{\rho} = \nu \nabla p(u)$.
\begin{proof}
 We have that
 \begin{align*}
 \mathcal{M}^\oplus(\alpha,\beta;y(u),p(u),\nu \nabla y(u),\nu \nabla p(u),u)
 = &\, \frac{\alpha}{2}\|y(u) - y_d\|^2,
\end{align*}
 which is zero if we let $\alpha$ go to zero.
\end{proof}
\end{proposition}
%--------------------------------------------------------------------------------------------------
 Although the majorant $\mathcal{M}^{\oplus}$ is a guaranteed,
 computable and sharp upper estimate for the
 discretization error in the combined norm, it only decreases with order $h$ when
 discretizing the mesh. However,
 the combined norm $|||\cdot|||$ is an $L^2$-norm, and, hence, decreases with order $h^2$.
 So the majorant $\mathcal{M}^{\oplus}$ is an overestimation for
 the combined norm.
Now, we introduce another norm which is a weighted
$H^1$-norm for the state depending 
on the corresponding control.
More precisely, we derive an estimate for the discretization error measured in the following norm:
\begin{align}
 \label{definition:combinedNormH1}
 |||u-v|||_1^2 := \frac{1}{2} \|y(u) - y(v)\|^2 + \frac{2 \lambda \underline{\nu}^2}{C_F^2}
 \|\nabla y(u) - \nabla y(v)\|^2.
\end{align}
%--------------------------------------------------------------------------------------------------
\begin{theorem}
 \label{theorem:DiffCombNormH1}
 For any control function $v \in U_{\text{ad}}$,
 we have the estimate 
 \begin{align}
  \label{equation:DiffCombNormH1}
  |||u-v|||_1^2 \leq 
  \mathcal{J}(y(v),v) - \mathcal{J}(y(u),u)
  + \frac{3 \lambda}{2 C_F^2} \left(\|\boldsymbol{\tau}-\nu \nabla \eta\|
 + C_F \|f+ v + \emph{div} \, \boldsymbol{\tau}\|\right)^2.
 \end{align}
\begin{proof}
Let $\delta > 0$ be an arbitrary but fixed parameter.
 Adding and subtracting
 $\nabla \eta$ as well as applying triangle inequality
 for $\frac{\underline{\nu}^2}{C_F^2 \delta}
 \|\nabla y(u) - \nabla y(v)\|^2$,
 we derive the following estimate:
 \begin{align*}
 \frac{\underline{\nu}^2}{C_F^2 \delta} \|\nabla y(u) - \nabla y(v)\|^2
 \leq \frac{\underline{\nu}^2}{2 C_F^2 \delta}
 \left(\|\nabla y(u) - \nabla \eta\|^2 + \|\nabla y(v) - \nabla \eta\|^2\right).
 \end{align*}
 Using (\ref{definition:majorantBVP}), adding and subtracting $v$ as well as
 applying twice triangle inequality, we arrive at the estimate
  \begin{align*}
 \frac{\underline{\nu}^2}{C_F^2 \delta} \|\nabla y(u) - \nabla y(v)\|^2
  \leq \frac{1}{2 C_F^2 \delta}
 \bigg(&\left(\|\boldsymbol{\tau}-\nu \nabla \eta\| 
 + C_F \|f+ v + \text{div} \, \boldsymbol{\tau}\| + C_F \|u -v \| \right)^2 \\
 &+ \left(\|\boldsymbol{\tau}-\nu \nabla \eta\|
 + C_F \|f+ v + \text{div} \, \boldsymbol{\tau}\|\right)^2\bigg) \\
   \leq \frac{3}{4 C_F^2 \delta}
 (\|&\boldsymbol{\tau}-\nu \nabla \eta\| 
 + C_F \|f+ v + \text{div} \, \boldsymbol{\tau}\|)^2
 + \frac{1}{4 \delta} \|u -v \|^2.
 \end{align*}
 By using (\ref{equation:DiffCombNorm}) and the previous estimate,
 we derive the inequality
   \begin{align*}
 |||u-v|||^2 &+ \frac{\underline{\nu}^2}{C_F^2 \delta} \|\nabla y(u) - \nabla y(v)\|^2
 - \frac{1}{4 \delta} \|u -v \|^2 \\
 &= \frac{1}{2} \|y(u) - y(v)\|^2
 + \frac{\underline{\nu}^2}{C_F^2 \delta} \|\nabla y(u) - \nabla y(v)\|^2 
 + \left(\frac{\lambda}{2} - \frac{1}{4 \delta} \right) \|u-v\|^2 \\
 &\leq \mathcal{J}(y(v),v) - \mathcal{J}(y(u),u) 
 + \frac{3}{4 C_F^2 \delta}
 \left(\|\boldsymbol{\tau}-\nu \nabla \eta\| 
 + C_F \|f+ v + \text{div} \, \boldsymbol{\tau}\|\right)^2.
 \end{align*} 
 Finally, choosing $\delta = 1/(2 \lambda)$
 yields
 the estimate (\ref{equation:DiffCombNormH1}).
\end{proof}
\end{theorem}
\begin{theorem}
\label{theorem:majorantCombNormH1}
For any $\zeta \in V$, let $v_{\zeta} \in U_{\text{ad}}$
be given by 
the projection formula (\ref{definition:projectionDiscrete}).
Then, we obtain the following error majorant:
\begin{align}
 \label{definition:majorantCombinedNormH1}
  |||u-v_\zeta|||_1^2 \leq \mathcal{M}^\oplus_1(\alpha,\beta;\eta,\zeta,\boldsymbol{\tau},\boldsymbol{\rho},v_{\zeta})
\end{align}
 for arbitrary
$\eta \in V$, 
$\boldsymbol{\tau}, \boldsymbol{\rho} \in H(\emph{div},\Omega)$
and 
$\alpha,\beta > 0$, where
\begin{align*}
 \mathcal{M}^\oplus_1(\alpha,\beta;\eta,\zeta&,\boldsymbol{\tau},\boldsymbol{\rho},v_{\zeta})
 = \frac{\alpha}{2}\|\eta - y_d\|^2
 + \frac{(1+\alpha)(1+\beta)C_F^2}{2\alpha\underline{\nu}^2} \|\boldsymbol{\tau}- \nu \nabla \eta\|^2 \\
 &+ \frac{(1+\alpha)(1+\beta)C_F^4}{2\alpha\beta\underline{\nu}^2} \|f+ v_{\zeta} + \emph{div} \, \boldsymbol{\tau}\|^2
 - \int_\Omega \left(f+ v_{\zeta} + \emph{div} \, \boldsymbol{\tau}\right) {\zeta} \, d\boldsymbol{x} \\
 &- \int_\Omega \left(\boldsymbol{\tau} - \nu \nabla \eta \right) \cdot \nabla {\zeta} \, d\boldsymbol{x}
 + \frac{1}{\underline{\nu}} \left(C_F \, \|\eta - y_d + \emph{div} \, \boldsymbol{\rho}\|
 + \|\boldsymbol{\rho} - \nu \nabla {\zeta}\| \right) \\
 &\times \left(\frac{C_F^3}{\lambda \underline{\nu}} \|\eta - y_d + \emph{div} \, \boldsymbol{\rho}\|
 + \frac{C_F^2}{\lambda \underline{\nu}} \|\boldsymbol{\rho} - \nu \nabla {\zeta}\| 
 + C_F \|f+ v_{\zeta} + \emph{div} \, \boldsymbol{\tau}\|
 + \|\boldsymbol{\tau} - \nu \nabla \eta \| \right) \\
 &+ \frac{3 \lambda}{2 C_F^2} \left(\|\boldsymbol{\tau}-\nu \nabla \eta\|
 + C_F \|f+ v_{\zeta} + \emph{div} \, \boldsymbol{\tau}\|\right)^2.
\end{align*}
\begin{proof}
 Applying (\ref{equation:DiffCombNormH1}) as well as
  (\ref{estimate:NotExactWithMajorant}) and (\ref{estimate:exactWithMinorant}) 
  finally leads to
  the estimate (\ref{definition:majorantCombinedNormH1}).
\end{proof}
\end{theorem}
%--------------------------------------------------------------------------------------------------
\begin{proposition}
\label{proposition:overallmajorantH1}
 The majorant
 $\mathcal{M}_1^\oplus(\alpha,\beta;\eta,\zeta,\boldsymbol{\tau},\boldsymbol{\rho},v_{\zeta})$
 defined in
 (\ref{definition:majorantCombinedNormH1}) attains the exact lower bound on the exact solution of
 the optimal control problem (\ref{equation:minfunc:OCP})-(\ref{equation:boxConstr:OCP}), or,
 equivalently, of
 the optimality system (\ref{equation:KKTSysVF:1})-(\ref{equation:KKTSysVF:3}),
 i.e.,
\begin{align*}
    \inf_{\substack{\eta,{\zeta} \in V, \boldsymbol{\tau},\boldsymbol{\rho} \in H(\emph{div},\Omega), \\
    v_{\zeta} = \mathbb{P}_{[u_a,u_b]} \{u_d - \frac{1}{\lambda} {\zeta} \}, \alpha,\beta > 0}}
  \mathcal{M}_1^\oplus(\alpha,\beta;\eta,\zeta,\boldsymbol{\tau},\boldsymbol{\rho},v_{\zeta}) = 0.
\end{align*}
 The infimum is attained for $v_\zeta=u$, $\eta = y(u)$, ${\zeta}= p(u)$,
 $\boldsymbol{\tau} = \nu \nabla y(u)$ and $\boldsymbol{\rho} = \nu \nabla p(u)$.
\begin{proof}
 We have that
 \begin{align*}
 \mathcal{M}_1^\oplus(\alpha,\beta;y(u),p(u),\nu \nabla y(u),\nu \nabla p(u),u)
 = &\, \frac{\alpha}{2}\|y(u) - y_d\|^2,
\end{align*}
 which is zero if we let $\alpha$ go to zero.
\end{proof}
\end{proposition}
%--------------------------------------------------------------------------------------------------

%%%%%%%%%%%%%%%%%%%%%%%%%%%%%%%%%%
\section{The Unconstrained Case}
\label{Sec5:UnconstrainedCase}
%%%%%%%%%%%%%%%%%%%%%%%%%%%%%%%%%%

In the unconstrained case, we have that $U_{\text{ad}} = L^2(\Omega)$, i.e.,
$v_\zeta$ satisfies the optimality condition
\begin{align}
 \label{equation:KKTSysVF:3UnconstrainedForv}
   \zeta + \lambda(v_\zeta - u_d) = 0 \qquad \qquad \text{in } \Omega,
\end{align}
or, equivalently, $v_\zeta = u_d - \frac{1}{\lambda} \zeta$ in $\Omega$.
The majorant (\ref{definition:majorant}) and
minorant (\ref{definition:minorant}) 
simplify to 
\begin{align}
\label{definition:majorantUnconstrained}
 \begin{aligned}
 \mathcal{J}^\oplus(\alpha,\beta;\eta,\zeta,\boldsymbol{\tau}):=
 &\frac{1+\alpha}{2}\|\eta - y_d\|^2
 + \frac{(1+\alpha)(1+\beta)C_F^2}{2\alpha\underline{\nu}^2} \|\boldsymbol{\tau}- \nu \nabla \eta\|^2 \\
 &+ \frac{(1+\alpha)(1+\beta)C_F^4}{2\alpha\beta\underline{\nu}^2} \|f+ u_d - \frac{1}{\lambda} \zeta + \text{div} \, \boldsymbol{\tau}\|^2
 + \frac{1}{2\lambda} \|\zeta\|^2,
 \end{aligned}
\end{align}
where $\mathcal{J}^\oplus(\alpha,\beta;\eta,\zeta,\boldsymbol{\tau})
 =  \mathcal{J}^\oplus(\alpha,\beta;\eta,\boldsymbol{\tau},v_\zeta)$,
and
\begin{align}
  \label{definition:minorantUnconstrained}
  \begin{aligned}
  \mathcal{J}^\ominus&(\eta,\zeta,\boldsymbol{\tau},\boldsymbol{\rho})
 = \frac{1}{2} \|\eta-y_d\|^2
 + \frac{1}{2\lambda} \|\zeta\|^2
 + \int_\Omega \left(f+ u_d - \frac{1}{\lambda} \zeta + \text{div} \, \boldsymbol{\tau}\right) {\zeta} \, d\boldsymbol{x} \\
 &+ \int_\Omega \left(\boldsymbol{\tau} - \nu \nabla \eta \right) \cdot \nabla {\zeta} \, d\boldsymbol{x} 
 - \frac{1}{\underline{\nu}} \left(C_F \, \|\eta - y_d + \text{div} \, \boldsymbol{\rho}\| + \|\boldsymbol{\rho} - \nu \nabla {\zeta}\| \right) \\
 &\times \left(\frac{C_F^3}{\lambda \underline{\nu}} \|\eta - y_d + \text{div} \, \boldsymbol{\rho}\|
 + \frac{C_F^2}{\lambda \underline{\nu}} \|\boldsymbol{\rho} - \nu \nabla {\zeta}\| 
 + C_F \|f+ u_d - \frac{1}{\lambda} \zeta + \text{div} \, \boldsymbol{\tau}\|
 + \|\boldsymbol{\tau} - \nu \nabla \eta \| \right),
  \end{aligned}
\end{align}
respectively.
We obtain the following estimate using minorant (\ref{definition:minorantUnconstrained}):
\begin{align}
\label{estimate:exactWithMinorantUnconstrained}
  \mathcal{J}(y(u),u) \geq \mathcal{J}^\ominus(\eta,\zeta,\boldsymbol{\tau},\boldsymbol{\rho})
  \qquad
  \forall \, \eta, \zeta \in V \,\, \forall \, \boldsymbol{\tau}, \boldsymbol{\rho} \in H(\text{div},\Omega).
\end{align}
The minorant $\mathcal{J}^\ominus$ is sharp, i.e.,
\begin{align}
 \label{definition:minorantSupUnconstrained}
  \sup_{\substack{\eta, \zeta \in V, \boldsymbol{\tau}, \boldsymbol{\rho} \in H(\text{div},\Omega)}}
  \mathcal{J}^\ominus(\eta,{\zeta},\boldsymbol{\tau},\boldsymbol{\rho})
  = \mathcal{J}(y(u),u),
\end{align}
 since the supremum of $\mathcal{J}^\ominus$
 is attained for the optimal state $\eta = y$ 
 and adjoint state ${\zeta}= p$,
 and their corresponding exact fluxes
 $\boldsymbol{\tau} = \nu \nabla y$
 and $\boldsymbol{\rho} = \nu \nabla p$.
 The optimal control is given by 
 $u = u_d - \frac{1}{\lambda} p$.
In the unconstrained case, the proof of Theorem~\ref{theorem:DiffCombNorm}
provides not only an inequality but even an equation. 
This result is presented
in the following theorem: 
\begin{theorem}
 For any control function $v \in L^2(\Omega)$, we have that
 \begin{align}
  \label{equation:DiffCombNormUnconstrained}
  |||u-v|||^2 = 
  \mathcal{J}(y(v),v) - \mathcal{J}(y(u),u).
 \end{align}
 \begin{proof}
 The equation (\ref{equation:DiffCombNormUnconstrained})
 follows by repeating the proof of Theorem~\ref{theorem:DiffCombNorm} with applying equation
 (\ref{equation:KKTSysVF:3Unconstrained}) instead of inequality (\ref{equation:KKTSysVF:3}).
 \end{proof}
\end{theorem}
The norm $|||\cdot|||$ defined in (\ref{definition:combinedNorm})
can be represented in terms of the state and the adjoint state (instead of the control) by
using the optimality condition (\ref{equation:KKTSysVF:3Unconstrained}) 
for $u$ and $v$ as follows:
\begin{align}
  \label{definition:combinedNormUnconstrained}
  \begin{aligned}
 |||u-v|||^2 &= \frac{1}{2} \|y(u) - y(v)\|^2 + \frac{\lambda}{2} \|u-v\|^2 \\
 &= \frac{1}{2} \|y(u) - y(v)\|^2 + \frac{\lambda}{2} \|u_d - \frac{1}{\lambda} p(u) - (u_d - \frac{1}{\lambda} p(v))\|^2 \\
 &= \frac{1}{2} \|y(u) - y(v)\|^2 + \frac{1}{2\lambda} \|p(u) - p(v)\|^2.
  \end{aligned}
\end{align}
Next, we present similar results to those of Theorem~\ref{theorem:majorantCombNorm}
and Proposition~\ref{proposition:overallmajorant} for the unconstrained case.
%--------------------------------------------------------------------------------------------------
\begin{theorem}
Let $U_{\text{ad}} = L^2(\Omega)$. For any $\eta, \zeta \in V$,
we obtain the following error majorant:
\begin{align}
 \label{definition:majorantCombinedNormUnconstrained}
  |||u-v_\zeta|||^2 \leq \mathcal{M}^\oplus(\alpha,\beta;\eta,\zeta,\boldsymbol{\tau},\boldsymbol{\rho})
  := \mathcal{J}^\oplus(\alpha,\beta;\eta,\zeta,\boldsymbol{\tau})
  - \mathcal{J}^\ominus(\eta,\zeta,\boldsymbol{\tau},\boldsymbol{\rho})
\end{align}
with $|||u-v_\zeta|||^2 = |||u-u_d + \frac{1}{\lambda} \zeta|||^2$ and
\begin{align*}
 \mathcal{M}^\oplus(\alpha,\beta&;\eta,\zeta,\boldsymbol{\tau},\boldsymbol{\rho})
 = \frac{\alpha}{2}\|\eta - y_d\|^2
 + \frac{(1+\alpha)(1+\beta)C_F^2}{2\alpha\underline{\nu}^2} \|\boldsymbol{\tau}- \nu \nabla \eta\|^2 \\
 &+ \frac{(1+\alpha)(1+\beta)C_F^4}{2\alpha\beta\underline{\nu}^2} \|f+ u_d - \frac{1}{\lambda} \zeta
 + \emph{div} \, \boldsymbol{\tau}\|^2
 - \int_\Omega \left(f+ u_d - \frac{1}{\lambda} \zeta + \emph{div} \, \boldsymbol{\tau}\right) {\zeta} \, d\boldsymbol{x} \\
 &- \int_\Omega \left(\boldsymbol{\tau} - \nu \nabla \eta \right) \cdot \nabla {\zeta} \, d\boldsymbol{x} 
 + \frac{1}{\underline{\nu}} \left(C_F \, \|\eta - y_d + \emph{div} \, \boldsymbol{\rho}\| + \|\boldsymbol{\rho} - \nu \nabla {\zeta}\| \right) \\
 &\times \left(\frac{C_F^3}{\lambda \underline{\nu}} \|\eta - y_d + \emph{div} \, \boldsymbol{\rho}\|
 + \frac{C_F^2}{\lambda \underline{\nu}} \|\boldsymbol{\rho} - \nu \nabla {\zeta}\| 
 + C_F \|f+ u_d - \frac{1}{\lambda} \zeta + \emph{div} \, \boldsymbol{\tau}\|
 + \|\boldsymbol{\tau} - \nu \nabla \eta \| \right),
\end{align*}
where $\boldsymbol{\tau},\boldsymbol{\rho} \in H(\emph{div},\Omega)$
are arbitrary and $\alpha,\beta > 0$.
\begin{proof}
 The estimate follows by applying (\ref{equation:DiffCombNormUnconstrained})
 and then using the estimates 
 (\ref{estimate:NotExactWithMajorant})
 and (\ref{estimate:exactWithMinorantUnconstrained}).
\end{proof}
\end{theorem}
%--------------------------------------------------------------------------------------------------
\begin{proposition}
\label{proposition:overallmajorantUnconstrained}
 Let $U_{\text{ad}} = L^2(\Omega)$.
 The majorant $\mathcal{M}^\oplus(\alpha,\beta;\eta,\zeta,\boldsymbol{\tau},\boldsymbol{\rho})$ defined in
 (\ref{definition:majorantCombinedNormUnconstrained}) attains the exact lower bound on the
 exact solution of
 the optimal control problem (\ref{equation:minfunc:OCP})-(\ref{equation:forwardpde:OCP}),
 or, equivalently, of
 the optimality system (\ref{equation:KKTSysVF:1}),(\ref{equation:KKTSysVF:2})
 and (\ref{equation:KKTSysVF:3Unconstrained}), i.e.,
\begin{align*}
  \inf_{\substack{\eta, \zeta \in V, \boldsymbol{\tau}, \boldsymbol{\rho} \in H(\emph{div},\Omega), \\
  \alpha,\beta > 0}} \mathcal{M}^\oplus(\alpha,\beta;\eta,\zeta,\boldsymbol{\tau},\boldsymbol{\rho}) = 0.
\end{align*}
 In other words, the infimum is attained for 
 the optimal state $\eta = y$ 
 and adjoint state $\zeta = p$
 with the optimal control $u = u_d - \frac{1}{\lambda} p$, 
 and for the exact fluxes $\boldsymbol{\tau} = \nu \nabla y$ and $\boldsymbol{\rho} = \nu \nabla p$.
\begin{proof}
 We have that
 \begin{align*}
 \mathcal{M}^\oplus(\alpha,\beta;y(u),p(u),\nu \nabla y(u),\nu \nabla p(u))
 = &\, \frac{\alpha}{2}\|y(u) - y_d\|^2,
\end{align*}
 which is zero if we let $\alpha$ go to zero.
\end{proof}
\end{proposition}
%--------------------------------------------------------------------------------------------------
Finally, we repeat the results of Theorems \ref{theorem:majorantCombNormH1}
and Proposition \ref{proposition:overallmajorantH1} for the unconstrained case.
%--------------------------------------------------------------------------------------------------
\begin{theorem}
\label{theorem:majorantCombNormH1Unconstrained}
Let $U_{\text{ad}} = L^2(\Omega)$.
For any $\eta,\zeta \in V$, we obtain the following error majorant:
\begin{align}
 \label{definition:majorantCombinedNormH1Unconstrained}
  |||u-v_\zeta|||_1^2 \leq \mathcal{M}^\oplus_1(\alpha,\beta;\eta,\zeta,\boldsymbol{\tau},\boldsymbol{\rho})
\end{align}
with $|||u-v_\zeta|||_1^2 = |||u-u_d + \frac{1}{\lambda} \zeta|||_1^2$ and
\begin{align*}
 \mathcal{M}^\oplus_1(\alpha,\beta&;\eta,\zeta,\boldsymbol{\tau},\boldsymbol{\rho})
 = \frac{\alpha}{2}\|\eta - y_d\|^2
 + \frac{(1+\alpha)(1+\beta)C_F^2}{2\alpha\underline{\nu}^2} \|\boldsymbol{\tau}- \nu \nabla \eta\|^2 \\
 &+ \frac{(1+\alpha)(1+\beta)C_F^4}{2\alpha\beta\underline{\nu}^2} \|f+ u_d - \frac{1}{\lambda} \zeta + \emph{div} \, \boldsymbol{\tau}\|^2
 - \int_\Omega \left(f+ u_d - \frac{1}{\lambda} \zeta + \emph{div} \, \boldsymbol{\tau}\right) {\zeta} \, d\boldsymbol{x} \\
 &- \int_\Omega \left(\boldsymbol{\tau} - \nu \nabla \eta \right) \cdot \nabla {\zeta} \, d\boldsymbol{x}
 + \frac{1}{\underline{\nu}} \left(C_F \, \|\eta - y_d + \emph{div} \, \boldsymbol{\rho}\|
 + \|\boldsymbol{\rho} - \nu \nabla {\zeta}\| \right) \\
 &\times \left(\frac{C_F^3}{\lambda \underline{\nu}} \|\eta - y_d + \emph{div} \, \boldsymbol{\rho}\|
 + \frac{C_F^2}{\lambda \underline{\nu}} \|\boldsymbol{\rho} - \nu \nabla {\zeta}\| 
 + C_F \|f+ u_d - \frac{1}{\lambda} \zeta + \emph{div} \, \boldsymbol{\tau}\|
 + \|\boldsymbol{\tau} - \nu \nabla \eta \| \right) \\
 &+ \frac{3 \lambda}{2 C_F^2} \left(\|\boldsymbol{\tau}-\nu \nabla \eta\|
 + C_F \|f+ u_d - \frac{1}{\lambda} \zeta + \emph{div} \, \boldsymbol{\tau}\|\right)^2
\end{align*}
for arbitrary
$\boldsymbol{\tau}, \boldsymbol{\rho} \in H(\emph{div},\Omega)$ 
and $\alpha,\beta > 0$.
\begin{proof}
 Applying (\ref{equation:DiffCombNormH1}) in the unconstrained case as well as
  (\ref{estimate:NotExactWithMajorant}) and (\ref{estimate:exactWithMinorantUnconstrained}) 
  finally leads to
  the estimate (\ref{definition:majorantCombinedNormH1Unconstrained}).
\end{proof}
\end{theorem}
%--------------------------------------------------------------------------------------------------
\begin{proposition}
\label{proposition:overallmajorantH1Unconstrained}
Let $U_{\text{ad}} = L^2(\Omega)$.
 The majorant
 $\mathcal{M}_1^\oplus(\alpha,\beta;\eta,\zeta,\boldsymbol{\tau},\boldsymbol{\rho})$
 defined in
 (\ref{definition:majorantCombinedNormH1Unconstrained}) attains the exact lower bound
 on the exact solution of
 the optimal control problem (\ref{equation:minfunc:OCP})-(\ref{equation:forwardpde:OCP}), or, equivalently, of
 the optimality system (\ref{equation:KKTSysVF:1}),(\ref{equation:KKTSysVF:2})
 and (\ref{equation:KKTSysVF:3Unconstrained}), i.e.,
\begin{align*}
    \inf_{\substack{\eta,{\zeta} \in V, \boldsymbol{\tau},\boldsymbol{\rho} \in H(\emph{div},\Omega), \\
    \alpha,\beta > 0}}
  \mathcal{M}_1^\oplus(\alpha,\beta;\eta,\zeta,\boldsymbol{\tau},\boldsymbol{\rho}) = 0.
\end{align*}
 The infimum is attained for
 the optimal state $\eta = y$ 
 and adjoint state $\zeta = p$
 with the optimal control $u = u_d - \frac{1}{\lambda} p$, 
 and for the exact fluxes $\boldsymbol{\tau} = \nu \nabla y$ and $\boldsymbol{\rho} = \nu \nabla p$.
\begin{proof}
 We have that
 \begin{align*}
 \mathcal{M}_1^\oplus(\alpha,\beta;y(u),p(u),\nu \nabla y(u),\nu \nabla p(u))
 = &\, \frac{\alpha}{2}\|y(u) - y_d\|^2,
\end{align*}
 which is zero if we let $\alpha$ go to zero.
\end{proof}
\end{proposition}

%%%%%%%%%%%%%%%%%%%%%%%%%%%%%%%%%%%%%%%%%%%%
\section{The Finite Element Discretization and the Preconditioned MINRES Solver}
\label{Sec6:FEDiscretization}
%%%%%%%%%%%%%%%%%%%%%%%%%%%%%%%%%%%%%%%%%%%%

In this section, we present the finite element discretization of the optimality system in
order to derive approximations of the state, control and adjoint state.
These approximations can be used for the computation of the majorants
and minorants. We start with the unconstrained case.

\subsection*{The unconstrained case}

Since the control can be eliminated from the optimality system
by using (\ref{equation:KKTSysVF:3Unconstrained}),
we only have to solve the system (\ref{equation:KKTSysVF:1})-(\ref{equation:KKTSysVF:2})
for the state $y \in V$ and the adjoint state $p \in V$.
For that, we approximate these unknown functions 
by finite element functions $y_h$, $p_h \in V_h \subset V$,
where the finite element space $V_h$ is defined as
\begin{align}
\label{def:Vh}
V_h = \mbox{span} \{\varphi_1, \dots, \varphi_n\}
\end{align}
with the standard nodal basis
$\{\varphi_i(\boldsymbol{x}) = \varphi_{ih}(\boldsymbol{x}): i=1,2,\dots,n_h \}$
and $h$ denotes the discretization parameter (mesh size)
such that $n = n_h = \mbox{dim} V_h = O(h^{-d})$.
Using continuous, piecewise linear finite elements on
triangles on a regular triangulation
to construct the finite element subspace $V_h$ and its basis
(see, e.g., \cite{WM:Ciarlet:1978}),
yields the following linear system:
\begin{align}
 \label{equation:DiscreteKKTUnconstrained}
 \left( \begin{array}{cc}
     M_h  &  -K_h \\
     -K_h  & -\lambda^{-1} M_h \end{array} \right) \left( \begin{array}{c}
     \underline{y}_h \\
     \underline{p}_h \end{array} \right) = \left( \begin{array}{c}
     \underline{y}_d \\
     -(\underline{f}+\underline{u}_d) \end{array} \right),
\end{align}
which has to be solved with respect to the nodal parameter vectors
$\underline{y}_{h} = (y_{i})_{i=1,\dots,n} \in \mathbb{R}^n$
and 
$\underline{p}_{h} = (p_{i})_{i=1,\dots,n} \in \mathbb{R}^n$
of the finite element approximations 
\begin{equation}
\label{equation:FourierCoeffFE}
y_{h}(\boldsymbol{x}) = \sum_{i=1}^n y_{i} \varphi_i(\boldsymbol{x}) \quad \mbox{and} \quad
p_{h}(\boldsymbol{x}) = \sum_{i=1}^n p_{i} \varphi_i(\boldsymbol{x})
\end{equation}
to the unknown functions $y(\boldsymbol{x})$ and $p(\boldsymbol{x})$.
The matrices $M_h$ and $K_h$ correspond to the mass and
stiffness matrices, respectively. Their entries are
computed by the following formulas:
\begin{align*}
\begin{aligned}
 M_h^{ij} = \int_{\Omega} \varphi_i \varphi_j \,d\boldsymbol{x}, \hspace{0.8cm}
 K_h^{ij} &= \int_{\Omega} \nu \, \nabla \varphi_i \cdot \nabla \varphi_j \,d\boldsymbol{x},
 \qquad \qquad i,j = 1,\dots,n.
\end{aligned}
\end{align*}
The right hand sides are given by
\begin{align*}
\begin{aligned}
 {\underline{y}_d} = \Big\lbrack \int_{\Omega} {y_d} \, \varphi_j \,d\boldsymbol{x} \Big\rbrack_{j=1,\dots,n}, \quad
 {\underline{u}_d} = \Big\lbrack \int_{\Omega} {u_d} \, \varphi_j \,d\boldsymbol{x} \Big\rbrack_{j=1,\dots,n}
 \quad \mbox{and} \quad
 \underline{f} = \Big\lbrack \int_{\Omega} f \, \varphi_j \,d\boldsymbol{x} \Big\rbrack_{j=1,\dots,n}.
\end{aligned}
\end{align*}
The discretized optimality system (\ref{equation:DiscreteKKTUnconstrained}) is a saddle point problem
and can be solved by a preconditioned MINRES (minimal residual) method, see
\cite{WM:PaigeSaunders:1975}.
Hence, it is crucial to construct preconditioners, which yield robust and fast
convergence for the preconditioned MINRES method.
In \cite{WM:SchoeberlZulehner:2007},
the following preconditioner was constructed:
\begin{align}
\label{definition:preconditionerPUnconstrained}
 \mathcal{P} = \left( \begin{array}{cc}
     M_h + \sqrt{\lambda} K_h & 0 \\
     0 & \frac{1}{\lambda}M_h + \frac{1}{\sqrt{\lambda}} K_h \end{array} \right)
\end{align}
leading to a robust convergence 
with respect to $h$ and 
$\lambda$ (as well as $\nu$ for our problem).

\subsection*{The constrained case}

In the case of having inequality constraints imposed on the control, we can reformulate the
variational inequality 
(\ref{equation:KKTSysVF:3}), which is equivalent to the projection formula
\begin{align*}
  u(\boldsymbol{x}) = \mathbb{P}_{[u_a(\boldsymbol{x}),u_b(\boldsymbol{x})]}
 \{u_d(\boldsymbol{x}) - \frac{1}{\lambda} p(\boldsymbol{x}) \},
\end{align*}
by introducing an additional 
parameter $\mu \in L^2(\Omega)$, 
see, e.g., \cite{WM:HinzePinnauUlbrichUlbrich:2009, WM:Troeltzsch:2010}.  
The idea is to apply
a primal-dual active set strategy in order to linearize the optimality system,
see, e.g., \cite{WM:HintermuellerItoKunisch:2002}.
The parameter $\mu$ is defined as $\mu = -\frac{1}{\lambda} p + u_d - u$ for
the exact solution $u \in U_{\text{ad}}$.
Moreover, if a function $v \in U_{\text{ad}}$ satisfies the relations
\begin{align*}
 v(\boldsymbol{x}) = \left\{
                \begin{array}{ll}
                  u_a(\boldsymbol{x}) & \text{ if } v(\boldsymbol{x}) + \mu(\boldsymbol{x}) < u_a(\boldsymbol{x}), \\
                  & \\
                  u_d(\boldsymbol{x}) - \frac{1}{\lambda} p(\boldsymbol{x}) 
                  & \text{ if } v(\boldsymbol{x}) + \mu(\boldsymbol{x}) \in
                  [u_a(\boldsymbol{x}),u_b(\boldsymbol{x})], \\
                  & \\
                  u_b(\boldsymbol{x}) & \text{ if } v(\boldsymbol{x}) + \mu(\boldsymbol{x}) > u_b(\boldsymbol{x}),
                \end{array}
             \right.
\end{align*}
then it is the optimal solution,
since this means that $v$ satisfies the projection formula.
In the following, we want to state the main steps of the primal-dual active set method
for the optimal control problem (\ref{equation:KKTSysVF:1})-(\ref{equation:KKTSysVF:3}).
Given an iterate 
$(y_{k-1},u_{k-1},p_{k-1},\mu_{k-1})$, the active and inactive sets are determined as follows:
\begin{align*}
 \mathcal{A}_{k}^a = \{\boldsymbol{x} \in \Omega: u_k(\boldsymbol{x}) + \mu_k(\boldsymbol{x})
 &< u_a(\boldsymbol{x}) \}, \quad 
 \mathcal{A}_{k}^b = \{\boldsymbol{x} \in \Omega: u_k(\boldsymbol{x}) + \mu_k(\boldsymbol{x})
 > u_b(\boldsymbol{x}) \}, \\
 &\mathcal{I}_{k} = \Omega \backslash (\mathcal{A}_k^a \cup \mathcal{A}_k^b).
\end{align*}
If $\mathcal{A}_{k}^a = \mathcal{A}_{k-1}^a$ and $\mathcal{A}_{k}^b = \mathcal{A}_{k-1}^b$,
then we have attained the optimal solution. Otherwise, the next iterate is the solution of the system
\begin{align*}
 - \text{div} (\nu(\boldsymbol{x}) \nabla y_{k}(\boldsymbol{x})) 
 &= f(\boldsymbol{x}) + u_{k}(\boldsymbol{x}), \,\,\,\, \boldsymbol{x} \in \Omega, \qquad
 y_{k}(\boldsymbol{x}) = 0, \,\, \boldsymbol{x} \in \Gamma, \\
  - \text{div} (\nu(\boldsymbol{x}) \nabla p_{k}(\boldsymbol{x})) 
 &= y_{k}(\boldsymbol{x}) - y_d(\boldsymbol{x}), \,\, \boldsymbol{x} \in \Omega, \qquad
 p_{k}(\boldsymbol{x}) = 0, \,\, \boldsymbol{x} \in \Gamma, \\
 u_{k}(\boldsymbol{x}) - u_d(\boldsymbol{x}) + \lambda^{-1} \chi_{\mathcal{I}_{k}} p_{k}(\boldsymbol{x})
 &= \chi_{\mathcal{A}_{k}^a} u_a(\boldsymbol{x}) + \chi_{\mathcal{A}_{k}^b} u_b(\boldsymbol{x}), \,\,
 \boldsymbol{x} \in \Omega,
\end{align*}
where $\chi$ denotes the characteristic function.
We eliminate again the control $u_{k}$ from the system
and obtain the following reduced optimality system written
in its variational formulation:
\begin{align*}
\begin{aligned}
   \int_\Omega y_{k} \, w \, d\boldsymbol{x}
   - \int_\Omega \nu \nabla p_{k} \cdot  \nabla w \, d\boldsymbol{x} 
   &= \int_\Omega y_d \, w \, d\boldsymbol{x}, \\
   \int_\Omega \nu \nabla y_{k} \cdot \nabla w \, d\boldsymbol{x}
   + \lambda^{-1} \int_{\mathcal{I}_{k}} p_{k} \, w \, d\boldsymbol{x}
   &= \int_\Omega (f + u_d) w \, d\boldsymbol{x} 
   + \int_{\mathcal{A}_{k}^a} u_a \, w \, d\boldsymbol{x}
   + \int_{\mathcal{A}_{k}^b} u_b \, w \, d\boldsymbol{x}, \\
\end{aligned}
\end{align*}
for all $w, q \in V$. Now, discretizing the reduced optimality system
by the finite element method with 
the space $V_h \subset V$ defined as in (\ref{def:Vh}) leads to the linear system
\begin{align}
 \label{equation:DiscreteKKT}
 \left( \begin{array}{cc}
     M_h  &  -K_h \\
     -K_h  & -\lambda^{-1} M_{h,\mathcal{I}_{k}} \end{array} \right) \left( \begin{array}{c}
     \underline{y}_{h,k} \\
     \underline{p}_{h,k} \end{array} \right) = \left( \begin{array}{c}
     \underline{y}_d \\
     -(\underline{f}+\underline{u}_d)
     - \underline{u}_a - \underline{u}_b \end{array} \right),
\end{align}
where the matrix $M_{h,\mathcal{I}_{k}}$ with the entries
\begin{align*}
\begin{aligned}
 M_{h,\mathcal{I}_{k}}^{ij} = \int_{\mathcal{I}_{k}} \varphi_i \varphi_j \,d\boldsymbol{x},
 \qquad \qquad i,j = 1,\dots,n,
\end{aligned}
\end{align*}
and the vectors
\begin{align*}
\begin{aligned}
 {\underline{u}_a} = \Big\lbrack \int_{\mathcal{A}_{k}^a} {u_a} \, \varphi_j \,d\boldsymbol{x} \Big\rbrack_{j=1,\dots,n}
 \quad \mbox{and} \quad
 {\underline{u}_b} = \Big\lbrack \int_{\mathcal{A}_{k}^b} {u_b} \, \varphi_j \,d\boldsymbol{x} \Big\rbrack_{j=1,\dots,n}
\end{aligned}
\end{align*}
have to be computed in every $k$-th iteration step of the level set method.
In order to solve the linear system (\ref{equation:DiscreteKKT})
by the preconditioned MINRES method, we need again a preconditioner,
which provides fast convergence and is as robust as possible.
Here, we refer to the preconditioner
\begin{align}
\label{definition:preconditionerP}
 \mathcal{P}_{\mathcal{I}_{k}} = \left( \begin{array}{cc}
     M_h + \sqrt{\lambda} K_h & 0 \\
     0 & \frac{1}{\lambda}M_{h,\mathcal{I}_{k}} + \frac{1}{\sqrt{\lambda}} K_h \end{array} \right),
\end{align}
which was presented in \cite{WM:Kollmann:2013}. For the preconditioner 
(\ref{definition:preconditionerP}), it can be proven robustness
with respect to the inactive set $\mathcal{I}_{k}$ and the mesh size $h$ (as well as
with respect to the diffusion parameter $\nu$ in our case), 
but not with respect to the cost parameter $\lambda$ (the condition number scales
like $1/\sqrt{\lambda}$).

% % %%%%%%%%%%%%%%%%%%%%%%%%%%%%%%%%%%%%%%%%%%%%
\section{Numerical Results}
\label{Sec7:NumericalResults}
% % %%%%%%%%%%%%%%%%%%%%%%%%%%%%%%%%%%%%%%%%%%%%

In this section, we present and discuss first numerical results for the
\textit{unconstrained case}.
The computational domain $\Omega = (0,1) \times (0,1)$ is uniformly
decomposed into triangles, and standard continuous, piecewise linear
finite elements are used for the discretization. In this case, the
Friedrichs constant is $C_F = 1/(\sqrt{2}\pi)$.

The construction of $\eta, \zeta, \boldsymbol{\tau}$ and
$\boldsymbol{\rho}$ is an important issue 
in order to obtain sharp 
guaranteed bounds from the majorants and minorants.
Since $\eta$ and $\zeta$ are constructed by continuous,
piecewise linear approximations,
their gradients are only piecewise constant.
Then, $\nabla \eta, \nabla \zeta \in L^2(\Omega)$, but
$\nabla \eta, \nabla \zeta \not\in H(\text{div},\Omega)$.
Hence, a flux reconstruction is needed in order to obtain
suitable fluxes $\boldsymbol{\tau}, \boldsymbol{\rho}
\in H(\text{div},\Omega)$.
A good reconstruction of the flux is an important and nontrivial topic.
We can regularize the fluxes by a post-processing operator
which maps the $L^2$-functions into $H(\text{div},\Omega)$,
see \cite{WM:Repin:2008}.
There are various techniques for realizing these post-processing steps such as, e.g.,
local post-processing by an elementwise averaging procedure
or by using Raviart-Thomas elements, see
\cite{WM:Repin:2008, WM:MaliNeittaanmaekiRepin:2014} and
references, therein.
In our numerical experiments, we use Raviart-Thomas elements of the lowest order,
see, e.g., \cite{WM:RaviartThomas:1977, WM:BrezziFortin:1991, WM:RobertsThomas:1991}.
We define the normal fluxes on interior edges $E_{mn}$ by
\begin{align*}
 (\boldsymbol{\tau} \cdot n_{E_{mn}})|_{E_{mn}} &= (\lambda_{mn} (\nabla \eta)|_{T_m}
 + (1-\lambda_{mn}) (\nabla \eta)|_{T_n})\cdot n_{E_{mn}},
\\
 (\boldsymbol{\rho} \cdot n_{E_{mn}})|_{E_{mn}} &= (\lambda_{mn} (\nabla \zeta)|_{T_m}
 + (1-\lambda_{mn}) (\nabla \zeta)|_{T_n})\cdot n_{E_{mn}},
\end{align*}
with $\lambda_{mn} = 1/2$ due to uniform discretization.
Here, $(\nabla \eta)|_{T_m}$, $(\nabla \zeta)|_{T_m}$, $(\nabla \eta)|_{T_n}$ and $(\nabla \zeta)|_{T_n}$
are constant vectors on two arbitrary, neighboring elements $T_m$ and $T_n$.
On boundary edges, the only one existing flux is used.
Hence, three normal fluxes are defined on the three sides of each element.
Inside, we reconstruct the fluxes $\boldsymbol{\tau}$ and $\boldsymbol{\rho}$
by the standard lowest-order Raviart-Thomas ($\text{RT}^0$-) extension of normal fluxes with
\begin{align*}
\text{RT}^0(\mathcal{T}_h) := \{&\boldsymbol{\tau} \in (L^2(T))^2:
\forall \, T \in \mathcal{T}_h \quad \exists \, a,b,c \in \mathbb{R}
\quad \forall \, \boldsymbol{x} \in T, \\
&\boldsymbol{\tau}(\boldsymbol{x}) = (a,b)^T + c \, \boldsymbol{x} \text{ and }
[\boldsymbol{\tau}]_E \cdot n_E = 0 \,\, \forall 
\text{ interior edges } E \},
\end{align*}
where $[\boldsymbol{\tau}]_E$ denotes the jump of $\boldsymbol{\tau}$ across the edge $E$
shared by two neighboring elements 
on a triangulation $\mathcal{T}_h$.
Altogether, the $\text{RT}^0$-extension yields averaged fluxes from $H(\text{div},\Omega)$, i.e.,
\begin{align*}
 \boldsymbol{\tau} = G_{\text{RT}}(\nabla \eta), \quad 
 \boldsymbol{\rho} = G_{\text{RT}}(\nabla \zeta), \quad
 G_{\text{RT}}:L^2(\Omega) \rightarrow H(\text{div},\Omega).
\end{align*}
Results for functional a posteriori error estimates
of parabolic time-periodic boundary value problems
using the same discretization and flux reconstruction techniques
can be found in \cite{WM:LangerRepinWolfmayr:2015}.

In order to solve the saddle point systems (\ref{equation:DiscreteKKTUnconstrained}),
we use the AMLI preconditioner proposed by Kraus in \cite{WM:Kraus:2012}
for an inexact realization of the block-diagonal preconditioner
(\ref{definition:preconditionerPUnconstrained})
in the MINRES method.
The numerical results where computed on grids of different mesh sizes (from $8 \times 8$
to $256 \times 256$).
The preconditioned MINRES iteration was stopped after
$8$ iteration steps in all computations using
the AMLI preconditioner
with $4$ inner iterations.
The presented CPU times in seconds $t^{\text{sec}}$ include the
computational times for computing the majorants and minorants,
which are very small in comparison to the
computational times of the solver.
All computations were performed on a laptop with
Intel(R) Core(TM) i5-4308U CPU @ 2.80GHz.

In the numerical experiment, we consider the following given data: 
the desired state
\begin{align*}
 y_d(\boldsymbol{x}) = (1 + 0.04 \pi^4) \sin(x_1 \pi) \sin(x_2 \pi),
\end{align*}
the desired control $u_d(\boldsymbol{x}) = 0$ and the source term $f(\boldsymbol{x}) = 0$.
We choose the parameter $\nu = 1$ and the cost parameter $\lambda = 0.01$.
The exact state and control are known and given by
\begin{align*}
 y(\boldsymbol{x}) = \sin(x_1 \pi) \sin(x_2 \pi) \qquad \text{and} \qquad
 u(\boldsymbol{x}) = 2 \pi^2 \sin(x_1 \pi) \sin(x_2 \pi),
\end{align*}
respectively. Hence,
the exact value of the cost functional is given by $\mathcal{J}(y(u),u) = 2.385$.
Table~\ref{tab1} presents 
the CPU times in seconds $t^{\text{sec}}$, the majorants (\ref{definition:majorantUnconstrained}) 
and minorants (\ref{definition:minorantUnconstrained}) as well as the efficiency indices
\begin{align*}
 I_{\text{eff}}^{\oplus} &=
 \frac{\mathcal{J}^{\oplus}(\alpha,\beta;\eta,\zeta,\boldsymbol{\tau})}{\mathcal{J}(y(u),u)},
\qquad 
\qquad \,\,\,\,
 I_{\text{eff}}^{\ominus} =
 \frac{\mathcal{J}(y(u),u)}{\mathcal{J}^{\ominus}(\eta,\zeta,\boldsymbol{\tau},\boldsymbol{\rho})}, \\
 I_{\text{eff}}^{\oplus / \ominus} &=
 \frac{\mathcal{J}^{\oplus}(\alpha,\beta;\eta,\zeta,\boldsymbol{\tau})}{\mathcal{J}^{\ominus}(\eta,\zeta,\boldsymbol{\tau},\boldsymbol{\rho})},
  \qquad \qquad
 I_{\text{eff}}^{\mathcal{M}_1} =
 \sqrt{\frac{\mathcal{M}_1^{\oplus}(\alpha,\beta;\eta,\zeta,\boldsymbol{\tau})}{|||u-v_{\zeta}|||_1^2}},
\end{align*}
obtained on grids of different mesh sizes.
The parameters $\alpha$ and $\beta$ are chosen by a proper minimization of the
majorant $\mathcal{J}^{\oplus}$
with respect to $\alpha$ and $\beta$, see \cite{WM:MaliNeittaanmaekiRepin:2014}.
In Table~\ref{tab2}, the individual parts of the majorants and minorants are presented.
They are denoted as follows:
\begin{align*}
\begin{aligned}
 \mathcal{R}_1^{\eta} &= \|\text{div} \, \boldsymbol{\tau} - \frac{1}{\lambda} \zeta\|, \quad 
 \mathcal{R}_1^{\zeta} = \|\eta - y_d + \text{div} \, \boldsymbol{\rho}\|, \quad
 \mathcal{R}_2^{\eta} = \|\boldsymbol{\tau}- \nabla \eta\|, \quad
 \mathcal{R}_2^{\zeta} = \|\boldsymbol{\rho}- \nabla \zeta\|, \\
 \mathcal{R}_3^{\eta} &= \|\eta - y_d\|, \,\,\,
 \mathcal{R}_3^{\zeta} = \|\zeta\|, \,\,\,
 \mathcal{R}_4 = \int_\Omega \left(\text{div} \, \boldsymbol{\tau}
 - \frac{1}{\lambda} \zeta \right) {\zeta} \, d\boldsymbol{x}, \,\,\,
 \mathcal{R}_5 =  \int_\Omega \left(\boldsymbol{\tau} - \nu \nabla \eta \right)
 \cdot \nabla {\zeta} \, d\boldsymbol{x}.
\end{aligned}
\end{align*}
\begin{table}[!ht]
\begin{center}
\begin{tabular}{|c|ccccccc|}
  \hline
   grid & $t^{\text{sec}}$ & $\mathcal{J}^{\ominus}$
   & $\mathcal{J}^{\oplus}$
   & $I_{\text{eff}}^{\oplus}$ 
   & $I_{\text{eff}}^{\ominus}$
   & $I_{\text{eff}}^{\oplus / \ominus}$
   & $I_{\text{eff}}^{\mathcal{M}_1}$  \\
  \hline
   $8     \times     8$ &  0.006  & 2.248 & 2.438 & 1.022 & 1.061 & 1.085 & 2.369 \\ %k=2
   $16   \times   16$ &  0.012  & 2.351 & 2.431 & 1.019 & 1.015 & 1.034 & 1.953 \\ %k=3
   $32   \times   32$ &  0.045  & 2.376 & 2.412 & 1.011 & 1.004 & 1.015 & 1.751 \\ %k=4
   $64   \times   64$ &  0.179  & 2.383 & 2.399 & 1.006 & 1.001 & 1.007 & 1.656 \\ %k=5
   $128 \times 128$ &  0.784  & 2.384 & 2.392 & 1.003 & 1.000 & 1.003 & 1.610 \\ %k=6
   $256 \times 256$ &  3.246  & 2.385 & 2.389 & 1.002 & 1.000 & 1.002 & 1.587 \\ %k=7
  \hline
\end{tabular}
\end{center}
\caption{Efficiency of the minorants and majorants.}
\label{tab1}
\end{table}
\begin{table}[!ht]
\begin{center}
\begin{tabular}{|c|cccccccc|}
  \hline
   grid & $\mathcal{R}_1^{\eta}$ & $\mathcal{R}_1^{\zeta}$ & $\mathcal{R}_2^{\eta}$ 
   & $\mathcal{R}_2^{\zeta}$ & $\mathcal{R}_3^{\eta}$ & $\mathcal{R}_3^{\zeta}$ & $\mathcal{R}_4$
   &  $\mathcal{R}_5$ \\
  \hline
   $8     \times     8$ & 1.146 & 0.235 & 0.053 & 0.011 & 1.928 & 0.094 & -0.013 & 0.001 \\ %k=2
   $16   \times   16$ & 0.606 & 0.121 & 0.016 & 0.003 & 1.943 & 0.098 & -0.004 & 0.000 \\ %k=3
   $32   \times   32$ & 0.307 & 0.061 & 0.004 & 0.001 & 1.947 & 0.098 & -0.001 & 0.000 \\ %k=4
   $64   \times   64$ & 0.154 & 0.030 & 0.001 & 0.000 & 1.948 & 0.099 & -0.000 & 0.000 \\ %k=5
   $128 \times 128$ & 0.077 & 0.015 & 0.000 & 0.000 & 1.948 & 0.099 & -0.000 & 0.000 \\ %k=6
   $256 \times 256$ & 0.039 & 0.008 & 0.000 & 0.000 & 1.948 & 0.099 & -0.000 & 0.000 \\ %k=7
  \hline
\end{tabular}
\end{center}
\caption{The indiviudal parts of the minorants and majorants.}
\label{tab2}
\end{table}

In both tables, we observe the efficiency of the AMLI preconditioned MINRES method presented in
\cite{WM:Kraus:2012, WM:LangerWolfmayr:2013}.
The computational times increase with a factor of four.
Moreover, one can see that $\mathcal{R}_1^{\eta}$ and $\mathcal{R}_1^{\zeta}$ reduce as a factor of two,
and $\mathcal{R}_2^{\eta}$  and $\mathcal{R}_2^{\zeta}$ as a factor of four
showing the efficiency of the applied flux reconstruction.
However, one could consider other 
flux reconstruction techniques in terms of improving 
the efficiency indices such as using higher order Raviart-Thomas elements, 
see \cite{WM:MaliNeittaanmaekiRepin:2014}.
Altogether we can observe that the majorants and minorants provide
good estimates for the value of the cost functional and can be used in order
to compute guaranteed estimates for the discretization error in the combined norms
as it is discussed and proved in this paper.

%%%%%%%%%%%%%%%%%%%%%%%%%%%%%%%%%%%%%%%%%%%%%%
\section{Conclusions}
\label{Sec8:Conclusions}
%%%%%%%%%%%%%%%%%%%%%%%%%%%%%%%%%%%%%%%%%%%%%%

This work was devoted to the derivation of guaranteed and fully
computable lower bounds (minorants)
for cost functionals of distributed elliptic optimal control problems
in order to close the gap of the already existing results on upper
bounds for these cost functionals, see \cite{WM:GaevskayaHoppeRepin:2006}.
An important result of this work was to prove 
that the discretization error in the state and the control can be computed
by the difference between majorant and minorant of the cost functional.
Altogether we derive a fully computable upper bound for the discretization error
in the state and the control,
which can be, in principle, used as object of direct minimization.
However, the aim of this work was not to present an algorithm for this minimization, but to derive  
first results on minorants for cost functionals of distributed elliptic optimal control problems
with control constraints as well as present first numerical tests for the theoretical
results derived.

\section*{Acknowledgments}
The author gratefully acknowledges the financial support by the
Austrian Academy of Sciences,
and thanks U. Langer, S. Repin and G. Wachsmuth for the fruitful discussions 
and valuable comments
regarding a posteriori error estimation. 

\bibliographystyle{siam}
\bibliography{WMPaper.bib}

\end{document}